\documentclass[journal,twoside,web]{ieeecolor}
\usepackage{generic}
\usepackage{cite}
\usepackage{amsmath,amssymb,amsfonts}
\usepackage{algorithmic}
\usepackage{graphicx}
\usepackage{algorithm,algorithmic}
\usepackage{hyperref}
\hypersetup{hidelinks}
\usepackage{textcomp}
\usepackage{ulem}
\def\BibTeX{{\rm B\kern-.05em{\sc i\kern-.025em b}\kern-.08em
    T\kern-.1667em\lower.7ex\hbox{E}\kern-.125emX}}
\usepackage{caption}
\usepackage{subcaption}

\newtheorem{theorem}{Theorem}[section]

\newtheorem{corollary}[theorem]{Corollary}
\newtheorem{definition}[theorem]{Definition}
\newtheorem{remark}[theorem]{Remark}
\newtheorem{assumption}[theorem]{Assumption}
\newtheorem{example}[theorem]{Example}

\newcommand{\R}{\mathbb{R}}

\newcommand{\sign}{\text{sgn}}

\begin{document}

\title{Neighboring Extremal Optimal Control Theory for Parameter-Dependent Closed-loop Laws}
\author{Ayush Rai, Shaoshuai Mou, and Brian D. O. Anderson
\thanks{This paragraph of the first footnote will contain the date on 
which you submitted your paper for review.}
\thanks{ Ayush Rai and Shaoshuai Mou are with the School of Aeronautics and Astronautics, Purdue University, West Lafayette, IN 47906 USA (e-mail: rai29@purdue.edu; mous@purdue.edu).}
\thanks{Brian D. O. Anderson is with the School of Engineering, Australian National University, Acton, ACT 2601, Australia (e-mail: Brian.Anderson@anu.edu.au).}}

\maketitle

\begin{abstract}
This study introduces an approach to obtain a neighboring extremal optimal control (NEOC) solution for a closed-loop optimal control problem, applicable to a wide array of nonlinear systems and not necessarily quadratic performance indices. The approach involves investigating the variation incurred in the functional form of a known closed-loop optimal control law due to small, known parameter variations in the system equations or the performance index. The NEOC solution can formally be obtained by solving a linear partial differential equation, akin to those encountered in the iterative solution of a nonlinear Hamilton-Jacobi equation. Motivated by numerical procedures for solving these latter equations, we also propose a numerical algorithm based on the Galerkin algorithm, leveraging the use of basis functions to solve the underlying Hamilton-Jacobi equation of the original optimal control problem. The proposed approach simplifies the NEOC problem by reducing it to the solution of a simple set of linear equations, thereby eliminating the need for a full re-solution of the adjusted optimal control problem. Furthermore,  the variation to the optimal performance index can be obtained as a function of both the system state and small changes in parameters,  allowing the determination of the adjustment to an optimal control law given a small adjustment of parameters in the system or the performance index. Moreover, in order to handle large known parameter perturbations, we propose a homotopic approach that breaks down the single calculation of NEOC into a finite set of multiple steps. Finally, the validity of the claims and theory is supported by theoretical analysis and numerical simulations.
\end{abstract}


%%%%%%%%%%%%%%%%%%%%%%%%%%%%%%%%%%%%%%%%%

\section{Introduction}
\label{sec:introduction}
\IEEEPARstart{O}{ptimal} control has garnered significant attention in recent decades due to its wide-ranging applications across fields such as aerospace, process control, economics, finance, and robotics. Its primary goal is to identify the most effective control inputs for a system characterized by specific dynamics, with the aim of optimizing a particular performance measure, often referred to as a performance index. Numerous methodologies have been developed to determine these optimal controls \cite{ anderson2007optimal, kirk2004optimal,lewis2012optimal,bryson2018applied}. Most of these approaches lack the ability to maintain optimality when faced with small changes in the parameters of the optimal control problem, other than by undertaking a complete reevaluation of the control inputs. Such parameter adjustments may arise from inherent problem requirements or the need for tuning to accommodate new or modified objectives\cite{jin2020pontryagin,jin2021safe,lu2022cooperative}. Our focus is on addressing this issue for closed-loop optimal control problems by seeking approaches to identify small adjustments in the control law that can provide close to optimal solutions even when confronted with small, known perturbations.

Neighboring extremal optimal control (NEOC) is the term referring to a systematic process for modifying an optimal control strategy to accommodate small perturbations in parameters. Such parameters may be in the system itself, or in the performance index, and originally were taken to be simply the initial or terminal conditions in the original optimal control problem. See in particular Breakwell et al.'s work in 1963, marking the initial contribution to the field of neighboring-extremal optimization techniques \cite{breakwell1963optimization}. The general problem is to develop a neighboring-optimal feedback control scheme, meaning an adjustment to an already known optimal control,  for an \textit{open-loop} optimal control problem in which the state variables are subject to initial and possibly terminal constraints. The prevailing method has been to use the second variation (linear-quadratic) theory to minimize the second variation of the performance index; it requires linearization of the system about a nominal trajectory and an associated `quadraticization' of the performance index. To this end, a Riccati transformation is often employed, along with a backward sweep method, to calculate linear feedback gains that can handle small state variations along the original optimal control path \cite{breakwell1963optimization,hymas1973neighboring, bryson2018applied}. 
In later work, \cite{lee1989neighbouring}, an extension of the above ideas was used to treat a dynamic optimization problem with variation in the parameters in the system equation or the performance index (as opposed to variable initial or terminal states). Again, the starting point is an open-loop solution of the underlying optimal control problem prior to parameter variation. This paper used a modified backward sweep method to obtain linear feedback laws defining the control variation as a function of the parameter variation, which  had to be small to secure the linearity of the laws. The problem context actually gave rise to the NEOC terminology.  

As noted, previous research in this field has largely focused on scenarios where the base optimal control problem generates an open-loop optimal control (even when the adjustment due to variation might be in closed-loop form). Closed-loop feedback however offers a significant advantage over open-loop control in dealing with the inevitable disturbances that can compromise the optimal performance of the system. In this work, we examine what might happen given a small parameter variation when the original optimal control problem has a solution in closed-loop form. More precisely, the solution is characterized using the Hamilton-Jacobi equation for the optimal performance, and an associated closed-loop feedback law for the optimal control, rather than an open-loop time function dependent on the initial state. For a large collection of systems and performance indices, specified in more detail later,  we aim to provide a theoretical tool for determining a change in the closed-loop control law resulting from a small perturbation in the parameters. This tool depends on the applicability of a precursor tool for iteratively solving the Hamilton-Jacobian for the base problem (a nonlinear partial differential equation of course) by using a sequence of linear partial differential equations. The change in the closed-loop control law gives a new law, termed the NEOC law, which is intended to approximate the exact law which would result from a solution of a new base problem involving the adjusted parameter value. 

The preliminary results of this work were presented in the conference paper \cite{rai23NEOC_CDC}. Compared to \cite{rai23NEOC_CDC}, this paper provides a rigorous theoretical analysis of the unperturbed and perturbed problems, and complex numerical simulations, and answers the following questions:
\begin{itemize}
    \item Does the recursive algorithm for unperturbed problems maintain the admissible behavior, whereby admissible behavior we essentially mean resulting in a finite performance index?
    \item What stability guarantee does the NEOC law provide?
    \item What is the error bound between the NEOC law versus optimal control law obtained after recalculation using a new base problem with adjusted parameter value?
    \item How does the approach simplify for an LQR system?
    \item How can one handle large parameter perturbations using the NEOC approach?
\end{itemize}
To the best of the authors' knowledge, this is the first investigation in this area and is applicable to a wide variety of problems, with commentary on some not explored here in the conclusions.

The paper is divided into 7 sections. Section \ref{sec:prob_formulation} aims to provide an overview of the problem formulation, by explaining the problem and high-level aspects of  the solution using a Hamilton-Jacobi equation and optimal control law formula associated with the original unperturbed problem. Because the solution to our problem of interest can be interpreted as a type of variation to a single iteration arising in one known particular approach to iteratively solving a Hamilton-Jacobi equation, the solution of the unperturbed problem and its solution via this approach is reviewed in Section \ref{sec:unpertubed} for the class of systems of interest in this paper. Following this,  Section \ref{sec:perturbed} covers the perturbed problem, with the previous section providing some kind of a template for the solution to the problem of interest. To this point, the solutions are all analytic or contained in formulas such as integrals over a semi-infinite time interval of trajectories of a known system. An actual numerical approach to solving the problem of interest using the Galerkin algorithm (which has earlier been widely used for Hamilton-Jacobi equation solution, see e.g. \cite{beard1997galerkin}) is detailed in Section \ref{sec:algorithm}, with illustrative examples provided in Section \ref{sec:examples}. Finally, the conclusions are provided in Section \ref{sec:conclusion}; this section also includes mention of future problems, including possibilities for adaptive control. 

 \textit{Notations:} We use the notation $\mathbb{R}^n$ and $\mathbb{R}^{n\times p}$ to denote the set of all $n \times 1$ real vectors and $n \times p$ real matrices, respectively. The transpose of a matrix or vector is denoted by $(\cdot)^{\top}$. The squared norm of a vector $v$ weighted with respect to a positive definite symmetric $R$ is denoted by $\| v\|^2_R$, and is equal to $v^{\top}Rv$. The gradient of a continuous differentiable scalar function $\bar h(\cdot)$ with respect to vector $x \in \mathbb{R}^n$ is defined as a column vector denoted by $\nabla_x \bar h(\cdot) \in \mathbb{R}^n$, while the derivative of a continuous differentiable vector function $h(\cdot)\in \mathbb{R}^n$,  with respect to vector $\alpha \in \mathbb{R}^q$ is defined using a Jacobian matrix $J_{h,\alpha}(\cdot) = [\partial h(\cdot)/\partial \alpha] \in \mathbb{R}^{n\times q}$ whose element in the $i$th row and $j$th column is $\partial h(\cdot)_i/\partial \alpha_j$. We use the notation $\{(\cdot)_l\}_{vec}$ to denote a column vector whose $l$th entry is given by $(\cdot)_l$, whereas we use $\{(\cdot)_l\}_{mat}$ to denote a matrix whose $l$th column is given by vector $(\cdot)_l$. The notation $B_r$ is used to define an open ball of radius $r$, i.e. $B_r = \{x\in \Omega | \|x\|<r\}$. Lastly, for two any functions $\omega_1, \omega_2,$ assumed square-integrable on a given set $\Omega$, we define the inner product as
\begin{equation} \label{eq:dot_prod}
\langle \omega_1,\omega_2\rangle_{\Omega}=\int_{\Omega}\omega_1(x)\omega_2(x)dx.
\end{equation}

\section{Problem formulation}
\label{sec:prob_formulation}

We set up a class of optimal control problems similar to those of \cite{moylan1973nonlinear}, but with a dependence on a parameter vector, call it $\alpha \in \mathbb{R}^{q}$. We assume that perturbations in  $\alpha$ are bounded. With this parameter dependence, the underlying control affine system dynamics is given as
\begin{equation}\label{eq:sys}
\dot x (t,\alpha)=f(x,\alpha)+g(x,\alpha)u \;\quad x(0,\alpha)=x_0(\alpha).
\end{equation}
where $x\in \Omega \subset \mathbb{R}^n$, $f(x,\alpha):\Omega \times \mathbb{R}^q \rightarrow \mathbb{R}^n$, $g(x,\alpha): \Omega \times \mathbb{R}^q \rightarrow \mathbb{R}^{n\times p}$, and $u(x,\alpha):\Omega \times \mathbb{R}^q \rightarrow \mathbb{R}^p$. We assume $f(\cdot)$ and $g(\cdot)$ are smooth, Lipschitz continuous on $\Omega$, and the equation has well-defined solutions in $\Omega$ for any closed-loop feedback control law $u(x,\alpha)$ of interest. We also assume $f(0,\alpha)=0$ for any $\alpha \in \mathbb{R}^q$ and that the system is completely controllable, in the sense that given any $(x(t_0,\alpha),t_0)$ and $(x(t_1,\alpha),t_1)$ there exists a smooth control $u$ defined on $[t_0,t_1]$ which will move the first state to the second. We assume that the set $\Omega$ is compact, and we restrict attention to trajectories (and associated control laws) which ensure it is invariant, i.e. $x(0,\alpha)\in\Omega$ means $x(t,\alpha)\in\Omega$ for all $t$. 

The performance index we consider, and we are interested in optimization for all initial conditions $x(0,\alpha)\in\Omega$ through a feedback law, is 

\vspace{-5pt}
\begin{small}
\begin{equation}\label{eq:perfindex}
\begin{split}
V(x_0,u(\cdot),\alpha)=\lim_{T\to\infty}\int_0^T[\|u(x(t,\alpha),\alpha)&\|_R^2+m(x(t,\alpha),\alpha)]dt \\ \quad&\mbox{s.t.}\quad    x(T,\alpha)=0,
\end{split}
\end{equation}
\end{small}where $m(x,\alpha)$ is a smooth Lipschitz, positive definite{\footnote{ Some relaxation of this assumption to permit some non-negative definite functions is possible, but for convenience in the subsequent analysis, we stay with the positive definiteness assumption.}}, radially increasing function that is strictly convex (has a positive definite Hessian) at the origin, while $R$ is a positive definite matrix. As shown in \cite[Chapter-3]{anderson2007optimal} and \cite{beard1997galerkin}, it is not enough for a system with a control law $u$ to be stabilizing for the limiting integral \eqref{eq:perfindex} to be finite. Therefore, it is necessary to introduce the concept of an \textit{admissible} control law.

\begin{definition} 
A control law $u: \Omega \times \mathbb{R}^{q} \rightarrow \mathbb{R}^{m}$ is considered admissible with respect to the performance index \eqref{eq:perfindex} for a given system dynamics \eqref{eq:sys} if it satisfies the following conditions: it is continuous on $\Omega$, $u(0,\alpha) = 0$, it stabilizes the system \eqref{eq:sys}, in the sense that $x(t)\to 0$ as $t\to\infty$, while also ensuring that $x(t,\alpha)\in\Omega\;\forall t$,  and it results in a finite integral in \eqref{eq:perfindex} for all $x(0,\alpha)$ in $\Omega$.
\end{definition}

We define the minimum performance index $\phi$ as the minimum of the cost function at optimal $u$ as\footnote{We assume the existence condition to hold true, meaning that there exists a control law $u^*$ that achieves the minimum of the performance index.}
\begin{equation}\label{eq:Phi}
    \phi(x,\alpha) = \min_u V(x,u(\cdot),\alpha).
\end{equation}
From the (steady-state) Hamilton-Jacobi equation \cite{anderson2007optimal}, for any given $\alpha$, we know that $\phi(x,\alpha)$ satisfies
\begin{equation}\label{eq:HJ}
\begin{split}
&[\nabla_x\phi(x,\alpha)]^{\top}f(x,\alpha)+m(x,\alpha) \\ -&\frac{1}{4}[\nabla_x\phi (x,\alpha)]^{\top}g(x,\alpha)R^{-1}g(x,\alpha)^{\top}\nabla_x\phi(x,\alpha)=0,
\end{split}
\end{equation}
and the optimal control law, which is provably stabilizing\footnote{By considering $\phi(x)$ as a Lyapunov function, one can demonstrate that the optimal control obtained from the Hamilton-Jacobi equation is stabilizing.} for the system \eqref{eq:sys}, i.e. ensures $x(t,\alpha)\to 0$ when $t\to\infty$, is given by
\begin{equation}\label{eq:oclaw}
u^*=-\frac{1}{2}R^{-1}g(x,\alpha)^{\top}\nabla_x \phi(x,\alpha).
\end{equation}

Suppose that the above calculations are done for a specific value of $\alpha$ (say $\alpha=\bar \alpha$), giving us unperturbed closed-loop feedback. Our task is to say what happens to the optimal performance index  $\phi$ and the control law \eqref{eq:oclaw} when $\alpha$ is changed to $\hat \alpha = \bar \alpha+\delta\alpha$ for some small known $\delta\alpha$. The problem can be regarded as one of establishing what the derivatives of $\phi$ and the optimal control law function are with respect to $\alpha$.

\section{Unperturbed problem}
\label{sec:unpertubed}
In this section, we begin by recalling one approach to solving the unperturbed problem where the parameter $\alpha$ is set to a fixed constant $\bar \alpha$. Since the parameter is fixed, we drop the $\alpha$ dependence for this section. Before delving into the solution, we provide a brief discussion of the previous literature surrounding the optimization problem involving \eqref{eq:sys} and \eqref{eq:perfindex} and its solution using \eqref{eq:HJ} and \eqref{eq:oclaw}. An early example of a simple first-order linear system with a quartic loss function in $x$, along with an explicit solution, is presented in \cite{bellman1964asymptotic}. For any system with scalar $x$, it is evident that \eqref{eq:HJ} defines the scalar quantity $\nabla_x \phi(x)$ up to a binary ambiguity, corresponding to the two solutions of a quadratic equation, and it is clear which one of the two achieves closed-loop stabilization. Additionally, \cite{bass1966optimal} demonstrates how an explicit solution can be easily obtained by carefully selecting a quartic loss function for a vector $x$. Although solving the equation is not straightforward, it has been extensively studied, as shown in \cite{markman2000iterative}.

Approximate solutions of the Hamilton-Jacobi equation can also be obtained, as suggested in \cite{leake1967construction}, by solving an iterative sequence of simpler equations, which in the limit yield the actual solution. (The calculation at a single iteration with the simpler equation serves to inspire our approach to solving the closed-loop NEOC problem.) In a broader sense, such approaches are known as policy iterations \cite{sutton2018reinforcement, howard1960dynamic}. The word `policy' refers to the control law in this case. The architecture of such algorithms involves policy evaluation and policy improvement. Starting with an admissible control input, in each iteration, one computes the associated performance index corresponding to the current control law (policy evaluation) and subsequently updates the control law based on the performance index (policy improvement). In more detail and for the particular form of system and performance index with which we are working, to start the iteration, one supposes the existence of an admissible (but not necessarily optimal) control $u_0=-\frac{1}{2}R^{-1}g^{\top}\nabla\phi_0(x)$ for some $\phi_0(x)$.   Subsequently, the iterative process ensures that at step $i$, a stabilizing control law, in the form of $u_i=-\frac{1}{2}g^{\top}\nabla\phi_i(x)$, is known and obtains the next iterate in a predefined manner. (The linear-quadratic version of this approach is actually known as the Kleinman algorithm, which replaces the solution of a time-invariant Riccati equation by the solution of a sequence of linear matrix equations,\cite{kleinman1968iterative}). More specifically, the algorithm defines $\phi_{i+1}(x)$ by the linear partial differential equation
\begin{equation}\label{eq:updatephi}
\begin{split}
[\nabla_x\phi_{i+1}(x)]^{\top}[f(x)-\frac{1}{2}g(x)R^{-1}g(x)^{\top}\nabla_x\phi_i(x)]= \\ -\frac{1}{4}[\nabla_x\phi_i(x)]^{\top}g(x)R^{-1}g(x)^{\top}\nabla_x\phi_i(x)-m(x)
\end{split}
\end{equation}

Because of the admissibility assumption, this is equivalent to setting
\begin{equation}\label{eq:phiformula}
\begin{split}
\phi_{i+1}(x)=\int_0^{\infty}\bigg[\frac{1}{4}[\nabla_y\phi_i(y)]^{\top}g(y)&R^{-1}g(y)^{\top}\nabla_y\phi_i(y) \\ &+m(y)\bigg]dt
\end{split}
\end{equation}
where the integration is performed along the trajectory with $y(\cdot)$ defined by
\begin{equation}\label{eq:y}
\dot y=f(y)-\frac{1}{2}g(y)R^{-1}g(y)^{\top}\nabla_y \phi_i(y)\quad y_i(0)=x.
\end{equation}

Evidently, $\phi_{i+1}(x)$ represents the value of the performance index when the closed-loop control $u_i(y)=-\frac{1}{2}R^{-1}g^{\top}(y)\nabla_y\phi_i(y)$ is used and the initial condition is $y(0)=x$.

Additional insight is obtained by considering the change from $\nabla_x\phi_i(x)$ to $\nabla_x\phi_{i+1}(x)$. For this purpose, define

\begin{small}
\begin{equation}\label{eq:HJerror}
\eta_i(x)= [\nabla_x\phi_i(x)]^{\top}f(x)+m(x) -\frac{1}{4}\big\|R^{-1}g(x)^{\top}\nabla_x\phi_i(x)\big\|_R^2
\end{equation}
\end{small}

The function $\eta_i$ can be interpreted as an error associated with $\nabla_x\phi_i$ being an approximate rather than exact solution of the steady state Hamilton-Jacobi equation. 

One straightforwardly obtains
\begin{align}
[\nabla_x\phi_{i+1}(x)-\nabla_x \phi_i(x)]^{\top}[f(x)-\frac{1}{2}g(x)&R^{-1}g(x)^{\top}\nabla_x\phi_i(x)] \nonumber\\& =-\eta_i(x)
\end{align}
which is equivalent to
\begin{equation}
\phi_{i+1}(x)-\phi_i(x)=\int_0^{\infty}\eta_i(y(s))ds
\end{equation}
with $y(\cdot)$ as in \eqref{eq:y} above.

The following theorem justifies the iteration process, being a particular instance of \cite{leake1967construction} and a generalization of the Kleinman algorithm\cite{kleinman1968iterative}.

\medskip 
\begin{theorem}
Consider an optimal control problem with system dynamics described by \eqref{eq:sys} and performance index defined in \eqref{eq:perfindex}, where $\alpha$ is a fixed constant. Given an admissible control input as an initial guess, the recursive algorithm outlined in \eqref{eq:updatephi} maintains the admissible behavior while also ensuring a monotonically decreasing sequence of performance indices $\phi_{i}(x), i=0.1.2\dots$.
\label{theorem1}
\end{theorem}

\noindent{\bf Proof:} 
This proof of this theorem is guided by the Kleinman lemma \cite{kleinman1968iterative} which we scale up for use in our nonlinear problem. We first verify the preservation of stability. Assume then that the control law $u_j = -\frac{1}{2}R^{-1}g^{\top} \nabla_x\phi_j(x)$ is stabilizing, which means, see \eqref{eq:phiformula} that $\phi_{j+1}(x)$ is a positive definite function; we will show that $u_{j+1} = -\frac{1}{2}R^{-1}g^{\top} \nabla_x\phi_{j+1}(x)$ is stabilizing. 
To simplify notation, we omit the dependence on $x$ for the functions $f$, $g$, $m$, and $\nabla_x\phi_j$, since they all depend only on $x$. Observe now that

\vspace{-9pt}

\begin{small}\begin{align}
\nabla&\phi_{j+1}^{\top}[f-\frac{1}{2}gR^{-1}g^{\top}\nabla\phi_{j+1}]\\ \nonumber
=&\nabla\phi_{j+1}[f-\frac{1}{2}gR^{-1}g^{\top}\nabla\phi_j] -\frac{1}{2}[\nabla\phi_{j+1}]^{\top}gR^{-1}g^{\top}\nabla\phi_{j+1} \\ \nonumber
  & +\frac{1}{2}\nabla\phi_{j+1}^{\top}gR^{-1}g^{\top}\nabla\phi_j\\\nonumber
=&-\frac{1}{4} [\nabla\phi_j^{\top}gR^{-1}g^{\top}\nabla\phi_j]-\frac{1}{2}[\nabla\phi_{j+1}]^{\top}gR^{-1}g^{\top}\nabla\phi_{j+1} \\  \nonumber
   & -m +\frac{1}{2}[\nabla\phi_{j+1}]^{\top}gR^{-1}g^{\top}\nabla\phi_j\\\nonumber 
=&- \frac{1}{4} [\nabla\phi_j-\nabla\phi_{j+1}]^{\top}gR^{-1}g^{\top}[\nabla\phi_j - \nabla\phi_{j+1}]\\ \nonumber
  &-\frac{1}{4}\nabla\phi_{j+1}^{\top}gR^{-1}g^{\top}\nabla\phi_{j+1}-m \nonumber
\end{align}
\end{small}
Recall that $\phi_{j+1}$ and $m$ are positive definite. Hence the above equation shows the asymptotic stability of $\dot x=[f-\frac{1}{2}gR^{-1}g^{\top}\nabla\phi_{j+1}]$ using $\phi_{j+1}$ as a Lyapunov function, i.e. the control law $u_{j+1}(x)$ is stabilizing, as required. Note that as $\phi_{j+1}$ is positive definite, it attains a minimum at the origin, and thus $\nabla_x\phi_{j+1}$ vanishes there, implying $u_{j+1}(0)=0$.

To show that the performance index $\phi_j(x)$ is monotonically decreasing with $j$, suppose that there are two stabilizing control laws given by
\begin{equation*}
  -\frac{1}{2}R^{-1}g^{\top}\nabla \phi_i \;\;\;\; \text{and}\;\;\;\; -\frac{1}{2}R^{-1}g^{\top}\nabla \phi_j.  
\end{equation*}

Next, observe that (using an expression similar to the left-hand side of \eqref{eq:updatephi})

\vspace{-10pt}

\begin{small}\begin{align}
[\nabla\phi_{j+1}]^{\top}[f-\frac{1}{2}&gR^{-1}g^{\top}\nabla\phi_i]\\ \nonumber
=[\nabla\phi_{j+1}]^{\top}&[f-\frac{1}{2}gR^{-1}g^{\top}\nabla\phi_j]\\ \nonumber
+&\frac{1}{2}[\nabla\phi_{j+1}]^{\top}gR^{-1}g^{\top}[\nabla\phi_j-\nabla\phi_i]\\\nonumber
=-\frac{1}{4}[\nabla\phi_j]^{\top}&gR^{-1}g^{\top}\nabla\phi_j-m \\ \nonumber +&\frac{1}{2}[\nabla\phi_{j+1}]^{\top}gR^{-1}g^{\top}[\nabla\phi_j-\nabla\phi_i]
\end{align}
\end{small}
Subtracting \eqref{eq:updatephi} then gives
\begin{small}
\begin{align}
[\nabla&\phi_{j+1}-\nabla\phi_{i+1}]^{\top}[f-\frac{1}{2}gR^{-1}g^{\top}\nabla\phi_i]\\ \nonumber
=&-\frac{1}{4}[\nabla\phi_j]^{\top}gR^{-1}g^{\top}\nabla\phi_j+\frac{1}{4}[\nabla\phi_i]^{\top}gR^{-1}g^{\top}\nabla\phi_i \\ \nonumber
& +\frac{1}{2}[\nabla\phi_{j+1}]^{\top}gR^{-1}g^{\top}[\nabla\phi_j-\nabla\phi_i]\\\nonumber
= & -\frac{1}{4}[\nabla\phi_j]^{\top}gR^{-1}g^{\top}\nabla\phi_j+\frac{1}{4}[\nabla\phi_i]^{\top}gR^{-1}g^{\top}\nabla\phi_i\\ \nonumber
&+\frac{1}{2}[\nabla\phi_{j+1}-\nabla\phi_i]^{\top}gR^{-1}g^{\top}[\nabla\phi_j-\nabla\phi_i] \\ \nonumber
&+\frac{1}{2}[\nabla\phi_i]^{\top}gR^{-1}g^{\top}[\nabla\phi_j-\nabla\phi_i]\\\nonumber
=&\frac{1}{2}[\nabla\phi_{j+1}-\nabla\phi_i]^{\top}gR^{-1}g^{\top}[\nabla\phi_j-\nabla\phi_i] \\\nonumber
&-\frac{1}{4}[\nabla\phi_j-\nabla\phi_i]^{\top}gR^{-1}g^{\top}[\nabla\phi_j-\nabla\phi_i]
\end{align}
\end{small}
Now suppose that we have first chosen the control law $-\frac{1}{2}R^{-1}g^{\top}\nabla \phi_j$ (corresponding to $j$). According to the algorithm presented, the next control law we choose is  $-\frac{1}{2}R^{-1}g^{\top}\nabla \phi_{j+1}$, where we have explained earlier how to go from $\phi_j$ to $\phi_{j+1}$. Let us identify in the immediately above formula the index $i$ with $j+1$. Then the equation gets simplified to
\begin{equation}
\begin{split}
[\nabla\phi_{j+1}-\nabla\phi_{j+2}]^{\top}[f-\frac{1}{2}gR^{-1}g^{\top}\nabla\phi_{j+1}]= \\-\frac{1}{4}[\nabla\phi_{j+1}-\nabla\phi_j]^{\top}gR^{-1}g^{\top}[\nabla\phi_{j+1}-\nabla\phi_{j}]
\end{split}
\end{equation}

Because the control law $-\frac{1}{2}R^{-1}g^{\top}\nabla \phi_{j+1}$ is stabilizing, this means that 
\vspace{-5pt}
\begin{equation}
\begin{split}
\phi_{j+1}-&\phi_{j+2} = \int_0^{\infty}\frac{1}{4}[\nabla\phi_{j+1}(y(s))-\nabla\phi_j(y(s))]^{\top} \\
 &\times gR^{-1}g^{\top}[\nabla\phi_{j+1}(y(s))-\nabla\phi_{j}(y(s))]ds
\end{split}
\end{equation}
with $y(\cdot)$ satisfying \eqref{eq:y}. Hence we have the crucial conclusion that $\phi_{j+2}(x) \leq \phi_{j+1}(x)\; \forall \; x$. This also means if $\phi_{j+1}(x)$ is finite $\forall \; x$, so will be $\phi_{j+2}(x)$ giving us admissibility of the associated control law. This completes the proof. $\hfill \blacksquare$

Of course, since the $\phi_i(x)$ are nonnegative for all $x$, the monotonicity property ensures the existence of a limiting function which is easily seen to satisfy the nonlinear Hamilton-Jacobi equation.

\section{Perturbation Problem}
\label{sec:perturbed}

In this section, we present the minor modifications required in the optimal control when small known perturbations affect the parameters of the optimal control problem. The results obtained in Section \ref{sec:unpertubed} give us a nominal closed-loop optimal feedback law, and they rely on repeatedly solving equations like \eqref{eq:updatephi}, which has a solution obtainable using \eqref{eq:phiformula} and \eqref{eq:y}. For the perturbation problem, a \textit{single} calculation parallel to \eqref{eq:updatephi}, \eqref{eq:phiformula}, and \eqref{eq:y} can be used. 
To study the consequence of small perturbation in the parameter, we define a vector function $\xi(x,\alpha) \in \mathbb{R}^{q}$ by
\begin{equation}
\xi(x,\alpha)= \nabla_\alpha\phi(x,\alpha),
\end{equation}
which also means that
\begin{equation}
J_{\xi,x}(x,\alpha) = J_{\nabla_\alpha \phi, x}(x,\alpha) .
\end{equation}
with the Jacobian matrix $J_{\xi,x}(x,\alpha) = [\partial \xi(x,\alpha)/ \partial x] \in \mathbb{R}^{q \times n}$.

Now simple differentiation of the parameterized steady-state Hamilton-Jacobi
equation \eqref{eq:HJ} yields
\begin{equation}
\begin{split}
&J_{\xi,x}(x,\alpha)\left[f(x,\alpha)-\frac{1}{2}g(x,\alpha)R^{-1}g(x,\alpha)^{\top}\nabla_x\phi(x,\alpha)\right] \\
&=- J_{f,\alpha}(x,\alpha)^{\top} \nabla_x\phi(x,\alpha) -\nabla_\alpha m(x,\alpha) \\ 
&+\frac{1}{2}\left\{\nabla_x\phi (x,\alpha)^{\top}\frac{\partial g(x,\alpha)}{\partial \alpha_l}R^{-1}g(x,\alpha)^{\top}\nabla_x\phi(x,\alpha)\right\}_{vec},
\end{split}
\label{eq:single_cal}
\end{equation}
where $\{(\cdot)\}_{vec}$ denotes a column vector whose $l$th entry is given by $(\cdot)\; \forall \; l = 1,\dots,q$. This means that formally there holds

\vspace{-10pt}

\begin{small}
\begin{equation}\label{eq:optPIderiv}
\begin{split}
&\xi(x,\alpha) = \int_0^{\infty} \Bigg(J_{f,\alpha}(y,\alpha)^{\top}\nabla_y\phi(y,\alpha) + \nabla_\alpha m(y,\alpha) \\ 
-&\frac{1}{2}\left\{ \nabla_y\phi (y,\alpha)^{\top}\frac{\partial g(y,\alpha)}{\partial \alpha_l}R^{-1}g(y,\alpha)^{\top}\nabla_y\phi(y,\alpha)\right\}_{vec} \Bigg) dt,
\end{split}
\end{equation}
\end{small}
with $y(\cdot)$ defined by
\begin{equation*}
\dot y=f(y,\alpha)-\frac{1}{2}g(y,\alpha)R^{-1}g(y,\alpha)^{\top}\nabla_y\phi(y,\alpha);\; y(0,\alpha)=x.
\end{equation*}

The change in optimum performance due to a small change $\delta \alpha$ away from an initial value $\bar \alpha$ is evidently given by $\xi(x,\bar \alpha)^{\top}\delta\alpha$ and the change in optimal control law is given by adding to the original feedback law the adjustment term
\begin{equation} 
\label{eq:control_adjust}
\begin{split}
    \delta u = -\frac{1}{2}R^{-1} &\Bigg(g(x,\bar\alpha)^{\top} J_{\xi,x}(x,\bar\alpha)^{\top} +  \\&\left\{\left[\frac{\partial  g(x,\bar\alpha)}{\partial \alpha_l}\right]^{\top} \nabla_x\phi(x,\bar\alpha)\right\}_{mat} \Bigg) \delta\alpha.
\end{split}
\end{equation} 
where $\{(\cdot)\}_{mat}$ denote a matrix whose $l$th column is given by $(\cdot)\; \forall \; l = 1,\dots,q$. Hence, the NEOC law is given as
\begin{align}
    u_{NE}(x,\hat \alpha) = u^{*}(x,\bar \alpha) + \delta u(x,\bar\alpha,\delta \alpha),
    \label{eq:control_NEOC}
\end{align}
where $\hat \alpha = \bar \alpha+\delta \alpha$ is the new parameter in the perturbed system, and $u^*$ is given by \eqref{eq:oclaw}.

For convenience, in our subsequent discussion, we will use the terms "original optimal control problem" and "new optimal control problem" to refer respectively to the unperturbed scenario with parameter $\bar \alpha$ and the scenario with parameter $\hat \alpha = \bar \alpha + \delta \alpha$. The NEOC law presented in \eqref{eq:control_NEOC} is designed to approximate the optimal control solution for this "new" problem. In essence, it provides the best possible estimate of the optimal control law for the new problem by considering the perturbation in parameters, denoted by $\delta \alpha$. One can conceptualize the closed-loop system of the new problem with NEOC law as a perturbation applied to the new system operating with its optimal control law. To illustrate this concept, let's consider the dynamics of the new optimal control problem, which can be expressed as:
\begin{align}
    \dot x (t,\hat \alpha)=f(x,\hat \alpha)+g(x,\hat \alpha)u. 
\end{align}
Consider $u^*(x,\hat \alpha)$ as the optimal control law for the perturbed problem, which is determined using \eqref{eq:oclaw}. When we apply the NEOC law, we obtain the following closed-loop dynamics
\begin{align}
    \dot x (t,\hat \alpha)=f(x,\hat \alpha)+g(x,\hat \alpha)u_{NE}, 
\end{align}
which can be formulated as the perturbed system of the new system as
\begin{align}
    \dot x (t,\hat \alpha)=f(x,\hat \alpha)&+g(x,\hat \alpha)u^*(\hat \alpha) \nonumber\\& + g(x,\hat \alpha)(u_{NE}-u^*(\hat \alpha)),
\end{align}
where $g(x,\hat \alpha)(u_{NE}-u^*(\hat \alpha))$ is the perturbation term, of course desirably small and vanishing as $\delta\alpha$ goes to zero. This enables us to discuss the stability guarantee of the NEOC law, as its stability directly relies on the stability guarantee of the optimal control law \eqref{eq:oclaw} \cite[Chapter-9]{khalil1992nonlinear}, which is determined by solving the Hamilton-Jacobi equation \eqref{eq:HJ}. Theorem \ref{theorem1} establishes asymptotic stability with the equilibrium at $x=0$. However, it is widely recognized that mere asymptotic stability (as opposed say to exponential asymptotic stability) is vulnerable to perturbations. Therefore, the question arises: Can we ensure more than just asymptotic stability for the optimal control law \eqref{eq:oclaw}? We provide a positive response to this question under certain assumptions by ensuring that, for any given $\alpha$, the optimal control law \eqref{eq:oclaw} guarantees uniformly globally asymptotic and locally exponential stability (UGALES) \cite{lin2019nonlinear}. To analyze local stability, it is necessary to linearize and quadraticize the problem around the origin, leading to a typical Linear Quadratic (LQ) regulator problem. This, in turn, demands adherence to the standard constraints \cite[Chapter 3]{anderson2007optimal} related to observability and controllability, which are assumed as follows.

\begin{assumption}
\label{assmp:lin_control}
For the linearized system of  \eqref{eq:sys} around the origin, the pair $\left[\frac{\partial f}{\partial x}(0,\alpha),g(0,\alpha)\right]$ is completely controllable.
\end{assumption}

\begin{assumption}
\label{assmp:lin_obs}  
For the linearized system of  \eqref{eq:sys} around the origin, the pair $\left[\frac{\partial f}{\partial x}(0,\alpha), Q^{1/2}(0,\alpha) \right]$ is completely observable, where $Q(0,\alpha)$ is the Hessian matrix of $m(x,\alpha)$ at the origin.\footnote{This assumption is trivially satisfied as $m$ is assumed to be strictly convex at the origin.}
\end{assumption}

\begin{theorem}
Consider an optimal control problem with system dynamics described by \eqref{eq:sys} and performance index defined in \eqref{eq:perfindex}, where $\alpha$ is a fixed constant. If the gradients of $f, g,$ and $m$ bounded and Lipschitz on $D = \{x\in \Omega \;| \; \|x\|<\rho\}$, then under Assumption \ref{assmp:lin_control} and Assumption \ref{assmp:lin_obs} the optimal control law \eqref{eq:oclaw} results in uniformly globally asymptotic and locally exponential stability of system \eqref{eq:sys} at the equilibrium $x=0$.
\label{theorem2}
\end{theorem}

\noindent{\bf Proof:} 
The UGALES property can be understood as the origin being exponentially stable in the local neighborhood of origin with uniform asymptotic stability in an arbitrarily large ball. We show that there exists a monotone non-decreasing continuous function $\beta(\cdot)>0$, and a real number $\lambda>0$, that satisfies
\begin{align}
    \|x(t,\alpha)\| \leq \beta(\|x_0(\alpha)\|)\|x_0(\alpha)\|e^{-\lambda t} \; \forall \; t \geq 0.
    \label{eq:ugales}
\end{align}

The closed-loop dynamics with optimal control law \eqref{eq:oclaw} is given as
\begin{align}
    \dot x (t,\alpha)=f(x, \alpha)-\frac{1}{2}g(x,\alpha)R^{-1}g(x,\alpha)^{\top}\nabla_x \phi(x,\alpha). 
\end{align}
The linearization of the dynamics about $x=0$ gives us
\begin{align*}
    \dot x (t,\alpha)=\left(\frac{\partial f(x, \alpha)}{\partial x}-\frac{1}{2}g(x,\alpha)R^{-1}g(x,\alpha)^{\top}J_{\nabla_{x} \phi(x,\alpha),x}\right)x.
\end{align*}
This occurs because $\nabla_{x} \phi(0,\alpha) = 0$ as $\phi(x,\alpha)$ is a positive definite and radially increasing function. The second derivative of $\phi(x,\alpha)$, $J_{\nabla_{x} \phi(x,\alpha),x}$, can be precisely calculated by differentiating the Hamilton-Jacobi equation \eqref{eq:HJ}. Under Assumption \ref{assmp:lin_control} and Assumption \ref{assmp:lin_obs} the linearized system is exponentially stable \cite{kalman1960contributions} due to the fact that $m$ is strictly convex at $x=0$. 

Following directly from Theorem 4.13 in \cite{khalil1992nonlinear} we can conclude that within a small ball $B_r$ around the origin (where $r$ is a small positive number\footnote{The specific value of $r$ depends on factors like $\rho$, system dynamics, and the performance index; please refer to Theorem 4.13 in \cite{khalil1992nonlinear}.}), the system exhibits exponential stability at $x=0$. This implies that, with certain positive constants $a$ and $\lambda$, we have the following relationship
\begin{align}
    \|x(t,\alpha)\| \leq a e^{-\lambda t}\|x_0(\alpha)\| \quad \forall\; t\geq 0, \|x(0,\alpha)\|\leq r. \label{eq:ugales_1}
\end{align}

From asymptotic stability of the system, for $\|x_0(\alpha)\| > r$, there is a first-time $T(x_0(\alpha))$ when the trajectory enters $B_r$. Such $T(x_0(\alpha)$ is finite, and because $\Omega$ is a compact set, $\max_{x_0(\alpha)\in \Omega}T(x_0(\alpha)$  is attained and finite. Then for $t \geq T$, and using local exponential stability we get 
\begin{align}
    \|x(t,\alpha)\| &\leq a e^{-\lambda (t-T)}\|x(T,\alpha)\|, \nonumber\\
    &\leq ar e^{-\lambda (t-T)} ,\nonumber\\
    &\leq a e^{-\lambda (t-T)}\|x_0(\alpha)\|. \label{eq:ugales_2}
\end{align} 

Lastly for $\|x_0(\alpha)\| > r$ and $t < T$, let $C(x(0,\alpha))=\max \{\|x(t,\alpha)\|   /  \|x(0,\alpha)\|\}$. Also, let $\bar C=\max_{x(0,\alpha)\in \Omega}C(x(0,\alpha))$ which is finite because $\Omega$ is compact. Then for $t < T$, we get
\begin{align}
    \|x(t,\alpha)\| &\leq \bar C \|x_0(\alpha)\| ,\nonumber\\
    &\leq \bar C \|x_0(\alpha)\| e^{-\lambda (t-T)}.
    \label{eq:ugales_3}
\end{align}
Combining \eqref{eq:ugales_1}, \eqref{eq:ugales_2}, and \eqref{eq:ugales_3} and setting
\begin{align*}
    \beta(\|x_0(\alpha)\|) = \max\{a,ae^{\lambda t}, \bar C e^{\lambda t}\},
\end{align*}
we get the desired result \eqref{eq:ugales}.

$\hfill \blacksquare$

\begin{remark}
    The outcome presented in Theorem \ref{theorem2} pertains to the stability of the nominal optimal control law at $x=0$. Considering that the closed-loop system using the NEOC law can be viewed as a perturbation of the new nominal optimal control system if the error between the NEOC and new optimal control law remains sufficiently small, it will maintain local exponential stability around the origin.
\end{remark}

Furthermore, our interest lies in establishing a quantitative measure of the approximation error between the NEOC law and the optimal control law for the new optimal control problem. In essence, we aim to determine the extent of perturbation that our approach can accommodate while maintaining a specified level of accuracy in the NEOC law. This accuracy can be gauged by assessing errors in various quantities of interest, such as the control law, performance index, trajectories, time to decay to zero, and more. The following theorem provides an error bound specifically for the control law. Although we don't present a formal proof for errors in other quantities, it is reasonable to assume that if the error in the control is minimal, errors in the performance index value and the trajectories themselves will also be minimal. This correlation is evident in our numerical simulations.

\begin{theorem}
\label{theorem3}
Consider an optimal control problem with system dynamics described by \eqref{eq:sys} and performance index defined in \eqref{eq:perfindex}. If $f, g, m$, and their gradients bounded and Lipschitz, then under the NEOC law obtained using \eqref{eq:control_NEOC}, the maximum permissible parameter perturbation, denoted by $\delta \alpha$, for a specified approximation error $\epsilon$ between the NEOC law and the optimal law of the new optimal control problem, $\|u_{NE}(x,\hat \alpha)- u^*(x,\hat \alpha)\| \leq \epsilon \; \forall \; x\in\Omega$, is given by
\begin{align*}
    \|\delta \alpha\| \leq \sqrt{\frac{2\epsilon}{M}},
\end{align*}
where $M$ is defined as
\begin{align*}
    M = \max_{x\in \Omega} \max_{\alpha\in [\bar \alpha, \hat \alpha]}\left\|\frac{\partial^2 u}{\partial \alpha^2}(x,\alpha)\right\|.
\end{align*}
\end{theorem}
\medskip

\noindent{\bf Proof:}
By employing the Taylor expansion for the optimal control law for the new problem, we observe that
\begin{align*}
    u^*(\bar \alpha+ \delta \alpha) = u^*(\alpha) + \frac{\partial u}{\partial \alpha}(\delta \alpha) + \frac{1}{2} \frac{\partial^2 u}{\partial \alpha^2}(\delta \alpha)^2+ \dots
\end{align*}
Breaking down the right-hand side of this expression, we find distinct components at play. Firstly, the initial term represents the nominal optimal control for the original problem. Secondly, the second term corresponds to the $\delta u$ term, which, in conjunction with the first term, constructs the NEOC solution. Lastly, the remaining higher-order terms encapsulate the approximation error characterizing the difference between the NEOC law and the optimal control law for the new optimal control problem. Designating the maximal magnitude of this approximation error concerning the parameter $\alpha$ as $R$, we can leverage the second derivative of the optimal control law and apply the Lagrange error bound to obtain
\begin{align*}
    R(x) \leq  \max_{\alpha\in [\bar \alpha, \hat \alpha]}\left|\frac{\partial^2 u}{\partial \alpha^2}(x,\alpha)\right|(\delta\alpha)^2.
\end{align*}
With $M = \max_{x\in \Omega} \max_{\alpha\in [\bar \alpha, \hat \alpha]}\left|\frac{\partial^2 u}{\partial \alpha^2}(x,\alpha)\right|$, for a given accuracy $\epsilon$, we get
\begin{align*}
    \frac{M(\delta \alpha)^2}{2} \leq \epsilon,
\end{align*}
resulting in the desired bound. Given that $f, g, m$ and their gradients are Lipschitz, it becomes apparent that within a compact set, the value of $M$ remains bounded.
$\hfill \blacksquare$

\section{Special Case: LQR system}
\label{sec:algorithm}
This section demonstrates the application of the NEOC approach to simple LQR (linear quadratic regulator) systems. We illustrate how the Riccati equation transforms into a set of linear equations. The system dynamics \eqref{eq:sys} can be simplified to a linear system as
\begin{align}
    \dot x= A(\alpha)x+ B(\alpha)u; \quad x(0,\alpha) = x_0(\alpha).
    \label{eq:lin_sys}
\end{align}
where $A(\alpha)\in \R^{n\times n} $ and $B(\alpha)\in \R^{n\times p} $. The performance index of interest can be expressed as a quadratic cost given by
\begin{align}
    V= \lim_{t\to \infty} \int_{0}^{T} [\|u(x(t,\alpha),\alpha)&\|_R^2+ \|x(t,\alpha)\|_{Q(\alpha)}^2]dt \nonumber\\ \quad&\mbox{s.t.}\quad    x(T,\alpha)=0,
    \label{eq:lin_cost}
\end{align}
where $Q$ and $R$ are positive definite matrices. In the case of the LQR system, for a given $\alpha$ the Hamilton-Jacobi equation simplifies to the well-known Riccati equation \cite{anderson2007optimal}, given by:
\begin{align}
A(\alpha)^{\top}P(\alpha)+ P(\alpha)A(\alpha) \nonumber -P(\alpha)&B(\alpha)R^{-1}B(\alpha)^{\top}P(\alpha) \\&+Q(\alpha)=0.
    \label{eq:ricati_equation}
\end{align}
and optimal control law can be given by linear feedback law
\begin{align}
    u^{*}(\alpha) = -R^{-1}B(\alpha)^{\top}P(\alpha)x.
    \label{eq:ricati_feedback}
\end{align}
Based on the developed perturbed problem in Section \ref{sec:perturbed}, under some small perturbation in the parameter $\alpha$ from some nominal value $\bar \alpha$, the change in the optimal control law can be obtained as
\begin{small}
\begin{align}
    \delta u = -R^{-1}\left(\left\{\frac{\partial B}{\partial \alpha_l}(\bar \alpha)^{\top}P(\bar \alpha)x + B(\bar \alpha)^{\top}\frac{\partial P}{\partial \alpha_l}(\bar \alpha)x\right\}_{mat}\right) \delta \alpha,
    \label{eq:ricati_NEOC}
\end{align}
\end{small}

\noindent where $\{(\cdot)\}_{mat}$ denote a matrix whose $l$th column is given by $(\cdot)\; \forall \; l = 1,\dots,q$. Here, we  have knowledge of $\frac{\partial B}{\partial \alpha_l}$, and we can calculate the quantity $\frac{\partial P}{\partial \alpha_l}$ by solving the following a set of linear equations\footnote{To simplify the notation, we have omitted the dependence on $\alpha$.}, obtained by taking the variation of the Riccati equation \eqref{eq:ricati_equation} with respect to $l^{th}$ element of $\alpha$:
\begin{align*}
     \left(A^{\top} -PBR^{-1}B^{\top}\right)\frac{\partial P}{\partial \alpha_l} 
     +\frac{\partial P}{\partial \alpha_l}\left(A-BR^{-1}B^{\top}P \right) 
     +\frac{\partial A}{\partial \alpha_l}^{\top}P\\ 
     + P\frac{\partial A}{\partial \alpha_l} 
     -P\frac{\partial B}{\partial \alpha_l}R^{-1}B^{\top}P -PBR^{-1}\frac{\partial B}{\partial \alpha_l}^{\top}P  + \frac{\partial Q}{\partial \alpha_l}=0,
\end{align*}
which can be simplified as
\begin{align}
    E^{\top} \frac{\partial P}{\partial \alpha_l} + \frac{\partial P}{\partial \alpha_l}E + F_l= 0
    \label{eq:lqr_lin_eq}
\end{align}
where $E = A-BR^{-1}B^{\top}P$ and $F_l = \frac{\partial A}{\partial \alpha_l}^{\top}P+ P\frac{\partial A}{\partial \alpha_l} -P\frac{\partial B}{\partial \alpha_l}R^{-1}B^{\top}P -PBR^{-1}\frac{\partial B}{\partial \alpha_l}^{\top}P  + \frac{\partial Q}{\partial \alpha_l}$. It is well known that $E$ is a Hurwitz matrix, hence these linear equations all have a unique solution. 

Solving the linear equation \eqref{eq:lqr_lin_eq} is generally a simpler task than solving the Riccati equation, which would be required if a new control law was computed starting with the new parameter value and not exploiting the calculation with $\bar\alpha$.

\section{Numerical Algorithm}
\label{sec:algorithm}
Although we have a NEOC solution in the form of equation \eqref{eq:optPIderiv} and \eqref{eq:control_adjust} that allows us to determine the sensitivity of the optimal performance index and optimal control law to a parameter change, obtaining an analytical or numerical solution can be challenging because it involves solving a linear partial differential equation. Nevertheless, given the similarity between the equations arising in our problem and those arising when iteratively solving the Hamilton-Jacobi equation, one can reasonably expect that one could apply to our problem a numerical approach commonly used for solving Hamilton-Jacobi equations. Successful policy iteration methods for optimal control problems have been implemented to solve Hamilton-Jacobi equations, including offline neural networks \cite{abu2005nearly}, sequential actor-critic neural networks \cite{vrabie2009adaptive}, and online actor-critic algorithms \cite{vamvoudakis2010online}. In \cite{beard1997galerkin}, researchers introduced an iterative algorithm based on the Galerkin spectral approximation method to solve a generalization of the Hamilton-Jacobi equation, viz. the generalized Hamilton-Jacobi-Bellman equation. This method offers a simple yet powerful approach to obtain the performance index in a functional form using basis functions. In this work, we extend and implement the Galerkin spectral approximation method to obtain the solution to the NEOC problem.

We explain first the original use of the method for an iterative Hamilton-Jacobi equation solution. It is assumed that in the iterative algorithm presented,  the $i$th approximant of the minimum performance index, denoted as $\phi_i(x,\alpha)$, can be represented as a sum of an infinite series of smoothly differentiable linearly independent weighted basis functions $\{\psi_j(x)\}_{j=1}^{\infty}$ along with their respective weights or coefficients $\{w_{ij}(\alpha)\}_{j=1}^{\infty}$, where the latter depends on the parameter $\alpha$. We select the basis functions such that $\phi_i(x,\alpha)$ is in the Hilbert space $L^2(\Omega)$, ensuring square integrable properties.\footnote{For more details refer \cite{beard1997galerkin}.} To make the computation feasible, we limit the summation of the infinite series to a finite number of terms which allows us to approximate $\phi_i(x,\alpha)$ to any desired degree of precision, giving us
\begin{equation}
\begin{split}
    \phi_i(x,\alpha) &\approx \sum_{j=1}^{N} w_{ij}(\alpha)\psi_j(x). \\
    & = w_i(\alpha)^{\top}\Psi(x),
\end{split}
\label{eq:phi_approx}
\end{equation}
where $w_i(\alpha) = [w_{i1}(\alpha), ..., w_{iN}(\alpha)]^{\top}\in \mathbb{R}^{N}$ and $\Psi(x) = [\psi_1(x),...,\psi_N(x)]^{\top} \in \mathbb{R}^{N}$. 

The associated optimal control law approximation is
\begin{equation}\label{eq:OCapprox}
u_i(x,\alpha)=-\frac{1}{2}R^{-1}g(x,\alpha)^{\top}J_{\Psi,x}(x)^{\top}w_i(\alpha)
\end{equation}
where $J_{\Psi,x}(x) = [\partial \Psi(x)/ \partial x] \in \mathbb{R}^{N\times n}$. 

The Galerkin projection method which we now describe constitutes a variant of the earlier iterative procedure of Section \ref{sec:unpertubed}; more precisely, assume that one has an expansion $w_i(\alpha)^{\top}\Psi(x)$ after $i$ steps as an approximation to the optimal performance index; one uses it to define the associated optimal control law approximation \eqref{eq:OCapprox}. Instead of then pursuing directly the calculation of $\phi_{i+1}(x,\alpha)$ as the associated performance index, one determines (with different calculations) an approximation $w_{i+1}(\alpha)^{\top}\Psi(x)$ to the index, which through appropriate choice of the coefficient vector $w_{i+1}(\alpha)$ is a least squares approximation over the whole set $\Omega$ to that index. Since $\phi_{i+1}(x,\alpha)$ in the normal iterative procedure is defined by

\begin{small}
\begin{equation*}
    \nabla_x\phi_{i+1}(x)^{\top}[f(x,\alpha)+g(x,\alpha)u_i(x,\alpha)]=-\|u_i(x,\alpha)\|^2_R-m(x,\alpha)
\end{equation*}
\end{small}
this means that $J_{\Psi,x}(x)^{\top}w_{i+1}(\alpha)$ must be an approximate solution of 
\begin{align}
    w_{i+1}(\alpha)^{\top}J_{\Psi,x}(x)[f(x,\alpha)+g(x,\alpha)u_i(x,\alpha)]\\\nonumber
    \approx-\|u_i(x,\alpha\|^2_R-m(x,\alpha)
\end{align}
The best least square approximate solution is obtained by choosing the coefficient vector $w_{i+1}(\alpha)$ to ensure that the error between the left side and the right side is orthogonal to the basis functions, i.e. for each $j=1,2,\dots,N$, there holds
\begin{align*}
    \Big\langle w_{i+1}(\alpha)^{\top}J_{\Psi,x}(x)[f(x,\alpha)+g(x,\alpha)u_i(x,\alpha)],\psi_j(x)\Big\rangle_{\Omega}\\\nonumber
    =-\langle\|u_i(x,\alpha\|^2_R+m(x,\alpha),\psi_j(x)\rangle_{\Omega},
\end{align*}
where the inner product is defined as per \eqref{eq:dot_prod}. There are $N$ linear equations in  $N$ unknowns, viz. the entries of $w_{i+1}(\alpha)$. It is established in \cite{beard1997galerkin} that the equation set is nonsingular, so that $w_{i+1}(\alpha)$ is well defined. 

The algorithm's convergence results as $N \rightarrow \infty$ and $i \rightarrow \infty$ are presented in \cite{beard1997galerkin}. In particular, one can assume that the sequence of iterations for some sufficiently large but fixed $N$ and some $\alpha=\bar \alpha$ will converge in practical terms once the iteration number $i$ reaches a value $k$. Then we can re-write the last equation as
\begin{equation}
\begin{split}
    &\left\langle w_{k}(\alpha)^{\top}J_{\Psi,x}(x)  f(x,\alpha), \psi_j(x) \right\rangle_\Omega + \langle  m(x,\alpha), \psi_j(x) \rangle_{\Omega} \\
    & +\frac{1}{4}\left\langle \left\|g(x,\alpha)^{\top}J_{\Psi,x}(x)^{\top}w_{k}(\alpha)\right\|^2_R, \psi_j(x) \right\rangle_\Omega  =0 ,
\end{split}   
\label{eq:galerkin_equality}
\end{equation}
for each $j=1,2,\dots,N$.

Moving now to the NEOC problem, these $N$ scalar equations will be the basis for determining the sensitivity of the optimal performance index and associated control law $($or more precisely, their approximations $w_k(\alpha)^{\top}\Psi(x)$ and $-\frac{1}{2}R^{-1}g(x,\alpha)^{\top}J_{\Psi,x}(x)^{\top}w_k(\alpha))$ to variation in the parameter $\alpha$.
In particular, by utilizing a calculation akin to that of \eqref{eq:single_cal}, we differentiate equation \eqref{eq:galerkin_equality} with respect to $\alpha \in \mathbb{R}^{q}$ to involve derivatives of the weighting coefficients with respect to the parameters. The following equation for each $j=1,2,\dots,N$, is to be understood as shorthand for $q$ scalar equations obtained by setting $ J_{(\cdot),\alpha_{\bar i}}$ in place of $J_{(\cdot),\alpha}$ for each $\bar i=1,2,\dots,q$.
The first term of \eqref{eq:galerkin_equality} can be differentiated as follows:
\begin{align*}
    &\left\langle J_{w_k,\alpha}(\alpha)^{\top} J_{\Psi,x}(x) f(x,\alpha), \psi_j(x) \right\rangle_\Omega \\&+ \left\langle J_{f,\alpha}(x,\alpha)^{\top}J_{\Psi,x}(x)^{\top} w_k(\alpha) , \psi_j(x) \right\rangle_\Omega.
\end{align*}
For the second term of \eqref{eq:galerkin_equality}, we have
\begin{align*}
    \left\langle \nabla_\alpha m(x,\alpha) , \psi_j(x) \right\rangle_\Omega.
\end{align*}
Lastly, the derivative of the third term of \eqref{eq:galerkin_equality} requires the formation of a tensor and can be simplified to
\begin{small}
\begin{align*}
    &-\frac{1}{2}\left\langle J_{w_k,\alpha}(\alpha)^{\top} J_{\Psi,x}(x) gR^{-1}g^{\top}J_{\Psi,x}(x)^{\top}w_k(\alpha),\psi_j \right\rangle_\Omega \\&-\frac{1}{2}\left\{\left\langle w_k(\alpha)^{\top} J_{\Psi,x}(x)\frac{\partial g}{\partial \alpha_l}R^{-1}g^{\top}J_{\Psi,x}(x)^{\top}w_k(\alpha),\psi_j \right\rangle_\Omega\right\}_{vec}.
\end{align*}
\end{small}
We observe that this differentiation of \eqref{eq:galerkin_equality} yields a linear equation set,  and it can be solved for    the entries of $J_{w_k,\alpha}=[\partial w_k(\alpha)/\partial \alpha] \in \mathbb{R}^{N\times q}$. These linear equations are provably solvable, as shown in Lemma 14 of \cite{beard1997galerkin}.

The variation in optimal control law can be obtained using \eqref{eq:control_adjust} under small perturbation $\delta\alpha$ from a constant value $\bar\alpha$, which yields
\begin{equation}
\label{eq:u_neoc_algo}
\begin{split}
    \delta u = -\frac{1}{2}R^{-1} \bigg(g(x,\bar\alpha)^{\top}J_{\Psi,x}(x)^{\top} J_{w_k,\alpha}(\bar\alpha) \\ +\left\{\left[\frac{\partial g(x,\bar\alpha)}{\partial \alpha_l}\right]^{\top} J_{\Psi,x}(x)^{\top} w_k(\alpha)  \right\}_{mat} \bigg) \delta\alpha.
\end{split}
\end{equation}

\section{Illustrative examples}
In this section, we present the NEOC approach results for four classical examples of closed-loop optimal control problems. For both the analytical and numerical approaches we compare three optimal control laws: the \textit{nominal} control law which is obtained for $\alpha = \bar\alpha$; the \textit{recalculated} optimal control law for the perturbed scenario where we re-run the complete algorithm using $\alpha = \bar \alpha + \delta\alpha$; and the \textit{NEOC} law where the nominal control law calculated for $\alpha = \bar \alpha$ is updated using a single calculation.

\label{sec:examples}
\begin{example}
    We first consider a simple LQR system, for which we use the linearized dynamics of an inverted pendulum mounted to on motorized cart with dynamics
    \begin{align*}
        \begin{bmatrix}
            \dot y \\ \ddot y \\ \dot \theta \\ \ddot \theta
        \end{bmatrix} = \begin{bmatrix}
            0 & 1 & 0 & 0\\ 
            0 & \frac{-(I+ml^2)b}{p} & \frac{m^2gl^2}{p} & 0\\ 
            0 & 0 & 0 & 1\\ 
            0 & \frac{-mlb}{p} & \frac{mgl(M+m)}{p} & 0
        \end{bmatrix}  \begin{bmatrix}
            y \\ \dot y \\ \theta \\ \dot \theta
        \end{bmatrix} + \begin{bmatrix}
            0 \\ \frac{I+ml^2}{p} \\ \dot \theta \\ \frac{ml}{p}
        \end{bmatrix}u,
    \end{align*}
\end{example}
\medskip
where $M=0.5$ and $m = 0.2$ are the masses of the cart and the pendulum, $b=0.1$ is the coefficient of friction for the cart, $l=0.3$ is the length of the pendulum, $I=0.006$ is the mass moment of inertia of the pendulum, and $y$ and $\theta$ are the cart position and pendulum angle from the vertically upward position. The value of $p = I(M+m)+Mml^2$. The performance index of this problem is of the form \eqref{eq:lin_cost}, with $Q= \begin{bmatrix}
    1 & 0 & 0 & 0\\
    0 & 0 & 0 & 0\\
    0 & 0 & 1 & 0\\
    0 & 0 & 0 & 0
\end{bmatrix}$, and $R=1$. We study the NEOC solution to this problem with perturbation $(\delta b=0.1)$ in the coefficient of friction for the cart. The optimal gain for the nominal control, recalculated control, and NEOC law (using \eqref{eq:ricati_NEOC}) are obtained as
\begin{align*}
    K_{nom} &= [-1.000, -1.656, 18.685, 3.459], \\
    K_{recal} &= [-1.000, -1.769, 18.732, 3.469], \\
    K_{neoc} &= [-1.000, -1.765, 18.716, 3.465]. 
\end{align*}
Evidently, the NEOC law, in which the adjustment to the initial control law is derived through the solution of a single linear equation, closely aligns with the recalculated control law.
\medskip

\begin{figure}[t!]
    \centering
    \begin{subfigure}[b]{0.4\textwidth}
        \centering        \includegraphics[width=\linewidth,trim={5cm 0.2cm 5.5cm 0.2cm},clip]{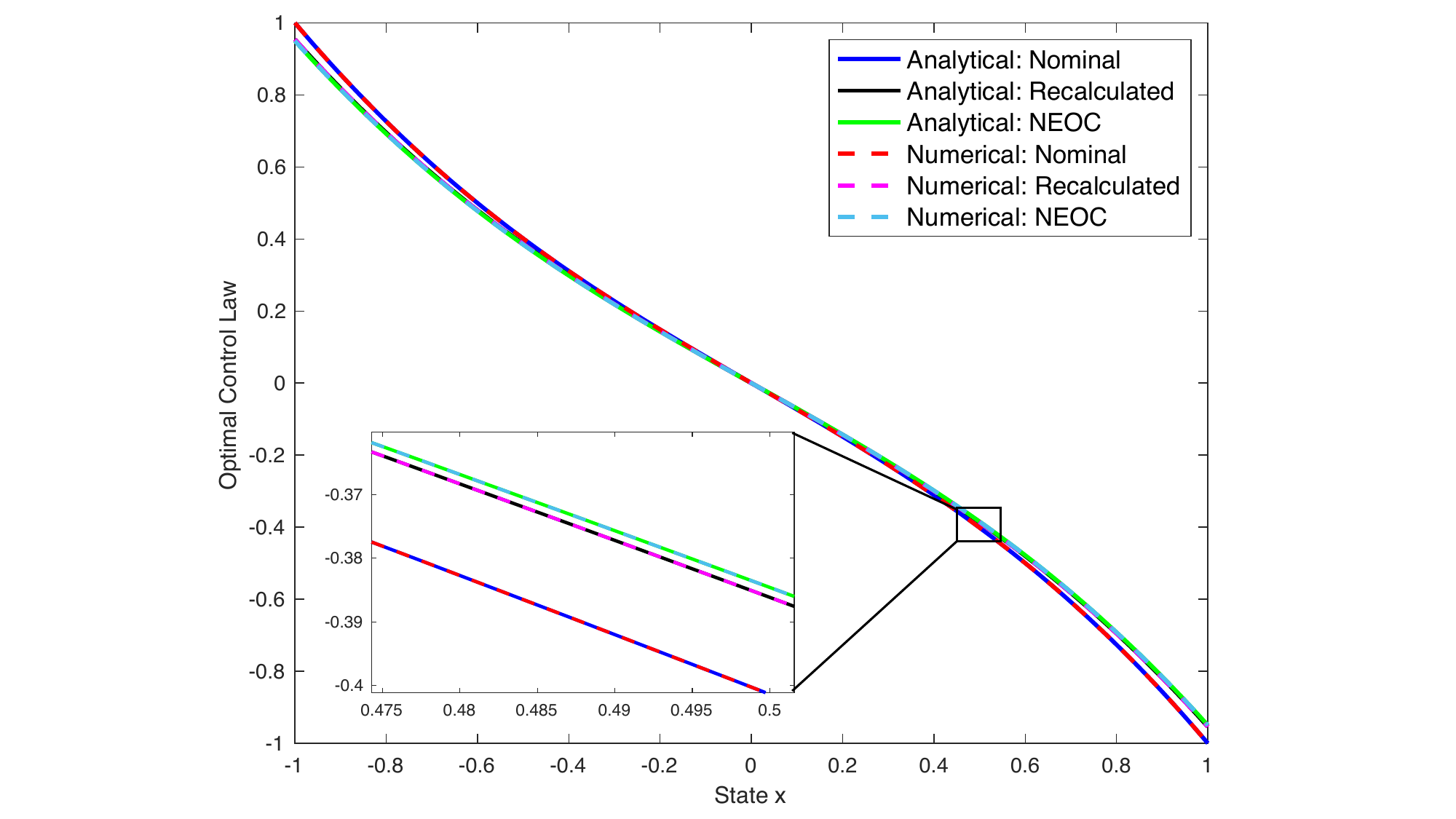} 
        \caption{NEOC performance} 
        \label{fig:example1_result}
    \end{subfigure}
    \hfill
    \begin{subfigure}[b]{0.4\textwidth}
        \centering
        \includegraphics[width=\linewidth,trim={0.8cm 0.5cm 1.6cm 0.5cm},clip]{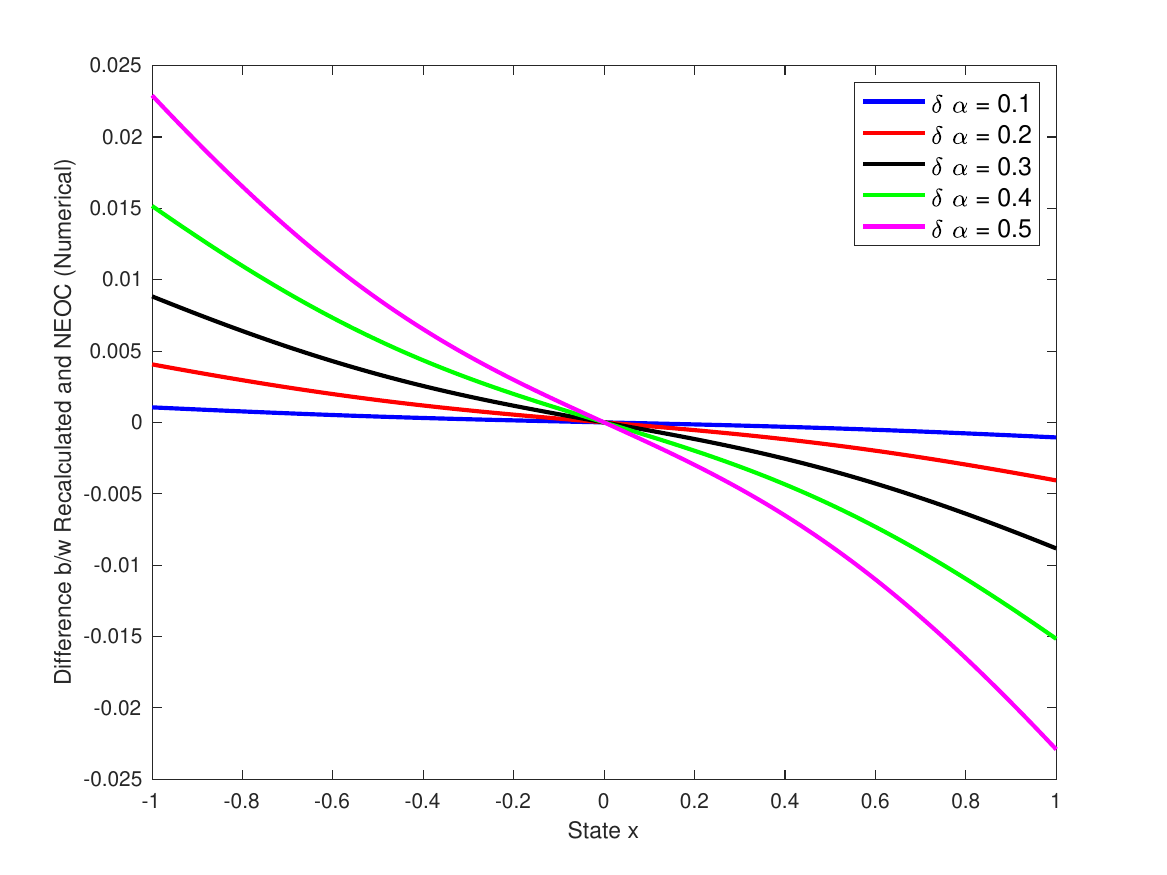}
        \caption{NEOC error}
        \label{fig:example1_error}
    \end{subfigure}
    \caption{Example 1- (a) Comparison of the analytical and numerical results for nominal, recalculated, and NEOC optimal control law for $\delta\alpha = 0.2$. (b) Difference between recalculated and NEOC laws for different values of perturbations.}
    \label{fig:example1}
\end{figure}

\begin{example}
For a simple non-linear system, we consider a general class of single-input single-output (i.e. scalar state) systems used in \cite{leake1967construction} with dynamics
\vspace{-5pt}
\begin{align*}
    \dot x(t,\alpha) = f(x,\alpha) + u,
\end{align*}
and a non-quadratic cost function-
\begin{small}
\begin{align*}
V=\lim_{T\to\infty}\int_0^T[u^2(x(t),\alpha)+m(x(t),\alpha)]dt \quad \mbox{s.t.}\quad x(T,\alpha)=0,
\end{align*}
\end{small}
where $m$ is a positive definite, radially increasing function that is strictly convex  at the origin.\footnote{Our results are concerned with closed-loop laws, so we can simply regard the initial state as fixed but arbitrary.}
\end{example}

In this case, the derivative $\nabla_x \phi(x)$ is a scalar, and the Hamilton Jacobi equation \eqref{eq:HJ} is simply a quadratic equation for $\nabla_x\phi(x)$. One of the two solutions is easily seen to define a stabilizing control law. In particular, we find that  
\begin{align*}
    \nabla_x\phi_{\infty}(x,\alpha) = 2f(x,\alpha)+2\sqrt{f(x,\alpha)^2+m(x,\alpha)},
\end{align*}
and the optimal control law is given using $u_{\infty}(x,\alpha) = -\frac{1}{2}\nabla_x\phi_{\infty}(x,\alpha)$ as
\begin{align*}
    u_{\infty}(x,\alpha) = -f(x,\alpha)-\sqrt{f(x,\alpha)^2+m(x,\alpha)}.
\end{align*}
This provides us with the analytical closed-loop optimal control law, which is expressible as a function of the parameters. In case of the presence of perturbation $\delta \alpha$ from the nominal value $\bar \alpha$, the optimal solution can be obtained by evaluating $u_{\infty}(x,\bar \alpha+ \delta \alpha)$. 

Using the NEOC approach, given that $u_{\infty}(x,\bar \alpha)$ has already been calculated, in the presence of known perturbation $\delta \alpha$ we add the adjustment term given by \eqref{eq:control_adjust}, which gives us
\begin{align*}
    \delta u_{\infty}(x,\bar\alpha,\delta \alpha) = -\bigg(\frac{\partial f}{\partial \alpha} +\frac{2\frac{\partial f}{\partial \alpha} +\frac{\partial m}{\partial \alpha}}{2\sqrt{f(x,\alpha)^2+m(x,\alpha)}}\bigg)\bigg|_{\alpha = \bar \alpha} \delta\alpha.
\end{align*}
This provides the analytical NEOC solution for the perturbation case as
\begin{align*}
    u_{NE}(x,\bar \alpha+\delta \alpha) = u_{\infty}(x,\bar \alpha) + \delta u_{\infty}(x,\bar\alpha,\delta \alpha).
\end{align*}

To assess the effectiveness of the numerical algorithm outlined in Section \ref{sec:algorithm}, we examine the performance of the algorithm for specific functions $f$ and $m$, where $f$ is defined as $-\alpha x$ and $m$ is defined as $(1+\alpha)x^2+ x^4$ and the value of $\alpha$ is centered at $1$. 

The corresponding analytical solution is given by
\begin{align*}
    u_{\infty}(x,\alpha) = \alpha x-x\sqrt{1+\alpha + \alpha^2+x^2},
\end{align*}
and the variation in the optimal control law around the nominal value of $\alpha = \bar \alpha$ can be expressed as
\begin{align*}
    \delta u_{\infty}(x,\bar\alpha,\delta \alpha) = \left(x-\frac{(1+2\bar \alpha)x}{2\sqrt{1+\bar \alpha+ \bar \alpha^2 + x^2}}\right)\delta \alpha.
\end{align*}

For the numerical algorithm, to approximate the minimum performance index, we utilize the basis function given by
\begin{align*}
    \{\psi_j\} = \{x^2 , x^4, x^6, x^8, x^{10}\}.
\end{align*}
To numerically obtain the nominal and recalculated control laws, we execute the Galerkin algorithm for 100 iterations with an initial admissible control input $u_0 = -5x$. The results obtained for the current examples are shown in Fig. \ref{fig:example1_result} for $\delta \alpha = 0.2$. The approximate optimal control laws obtained using NEOC are found to closely match the recalculated optimal control laws for both analytical and numerical cases. To further evaluate the accuracy of the NEOC method, we analyze the difference between the recalculated optimal and NEOC laws for different values of perturbations. As shown in Fig. \ref{fig:example1_error}, even after changing the parameter by $50\%$, the maximum NEOC error observed is just $0.023$. In this case, the analytical and numerical results cannot be visually distinguished because we are dealing with a simple example. However, in general, the accuracy of the results depends on the number of basis functions and iterations used to approximate $\phi$.

\begin{figure}[!htb]
    \centering
    \begin{subfigure}[b]{0.4\textwidth}
        \centering
        \includegraphics[width=\linewidth,trim={4.5cm 0.2cm 5cm 0.2cm},clip]{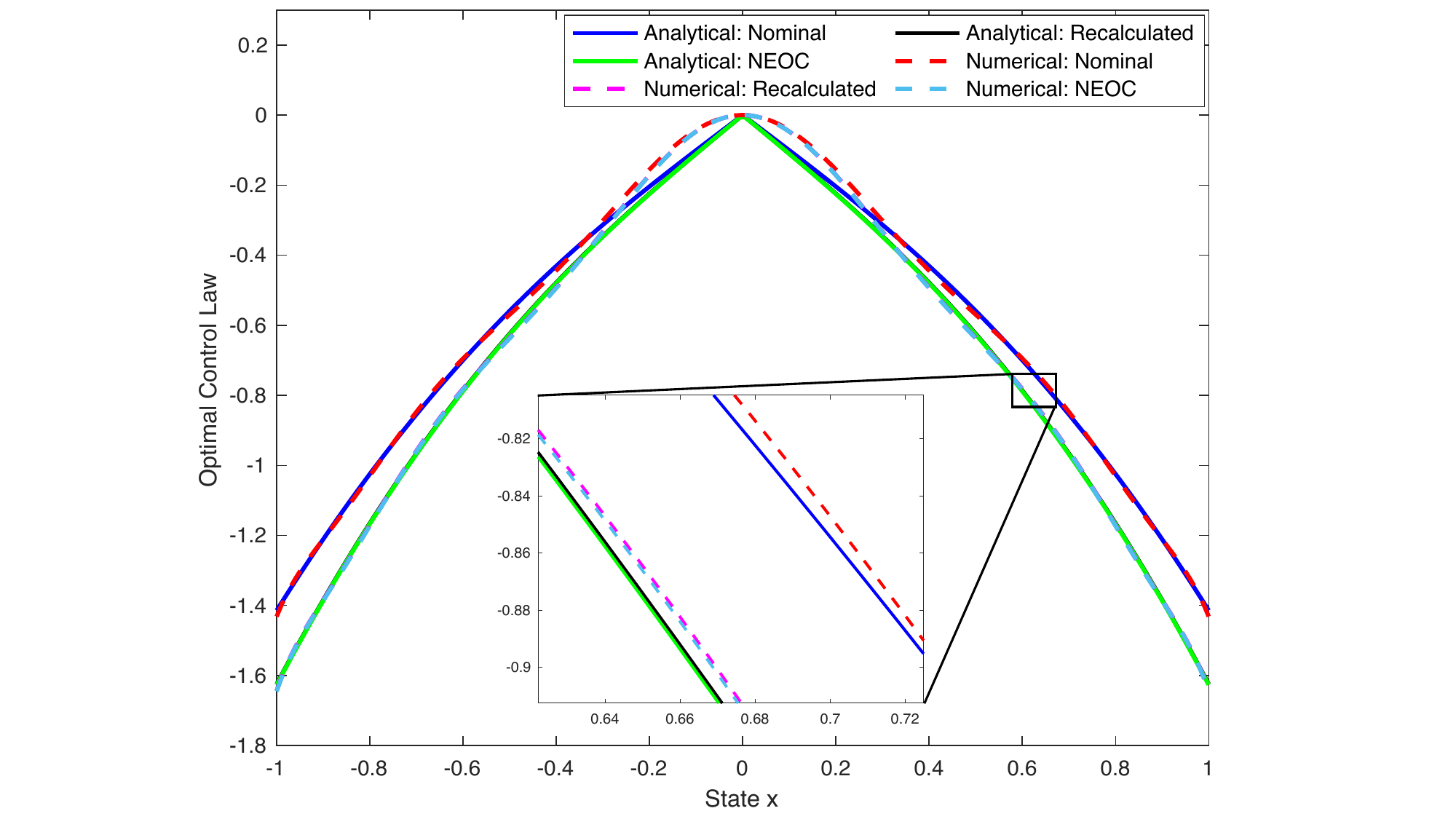}
        \caption{NEOC performace}
        \label{fig:example2_result}
    \end{subfigure}
    \hfill
    \begin{subfigure}[b]{0.4\textwidth}
        \centering
        \includegraphics[width=\textwidth,trim={0.9cm 0.5cm 1.2cm 0cm},clip]{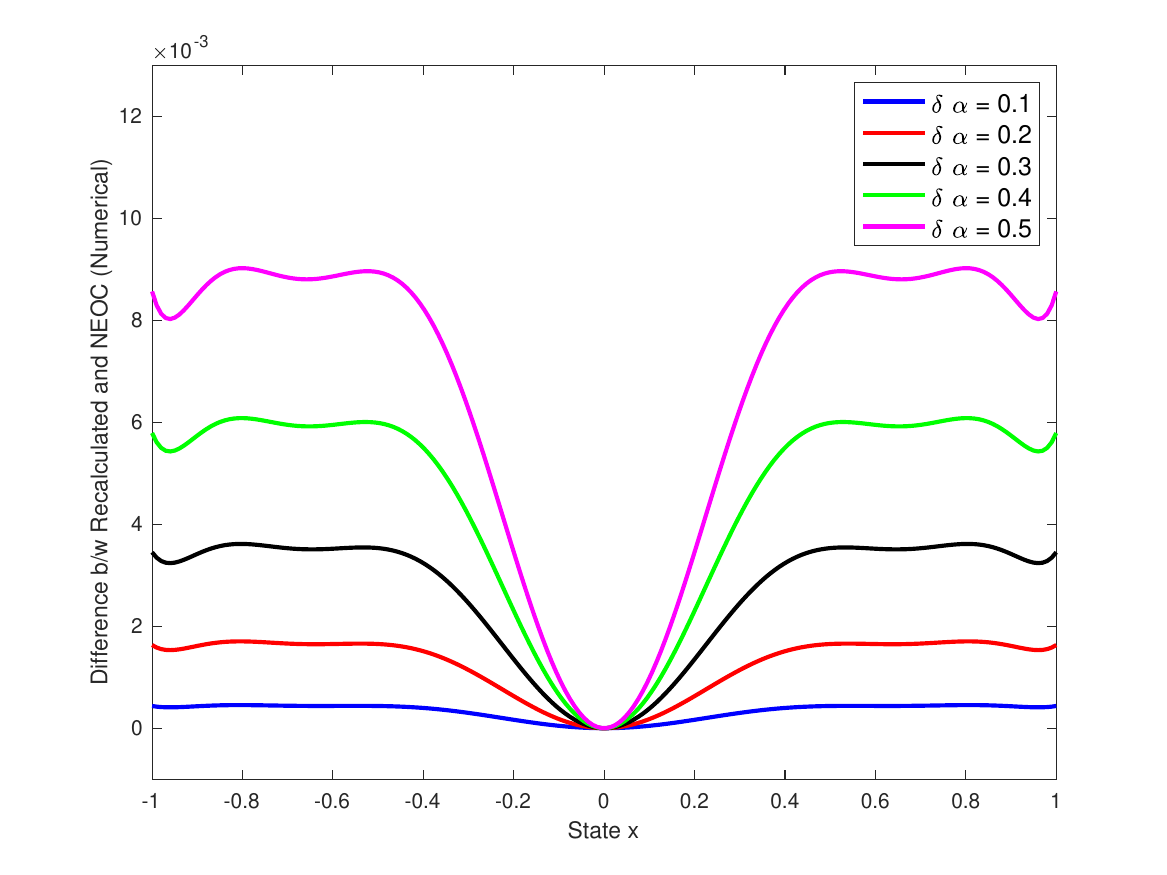}
        \caption{NEOC error}
        \label{fig:example2_error}
    \end{subfigure}
    \caption{Example 2- (a) Comparison of the analytical and numerical results for nominal, recalculated, and NEOC optimal control law for $\delta\alpha = 0.2$. (b) Difference between recalculated and NEOC laws for different values of perturbations.}
    \label{fig:example2}
\end{figure}

\begin{example}
We next consider a simplified class of bilinear systems used in \cite{beard1995improving}, given by
\begin{align*}
    \dot x(t,\alpha) = g(x,\alpha)u
\end{align*}
and a non-quadratic cost function-
\begin{small}
\begin{align*}
V=\lim_{T\to\infty}\int_0^T[u^2(x(t),\alpha)+m(x(t),\alpha)]dt \quad \mbox{s.t.}\quad x(T,\alpha)=0,
\end{align*}
\end{small}
where $m$ is a positive definite, radially increasing function that is strictly convex at the origin.
\end{example}

Similarly to the previous case, the Hamilton-Jacobi equation can be analytically solved, leading to the optimal solution given by 
\vspace{-2pt}
\begin{align*}
    \nabla_x\phi_{\infty}(x,\alpha) = \frac{ 2\sqrt{m(x,\alpha)}}{|g(x,\alpha)|}.
\end{align*}
This can be deduced utilizing equation \eqref{eq:updatephi}.
The optimal control law is given using $u_{\infty}(x,\alpha) = -\frac{1}{2}g(x,\alpha)\nabla_x\phi_{\infty}(x,\alpha)$ as-
\begin{align*}
    u_{\infty}(x,\alpha) = -\sign(g(x,\alpha))\sqrt{m(x,\alpha)},
\end{align*} 
providing us with the analytical closed-loop optimal law.

Using the NEOC approach, in the presence of a known perturbation $\delta \alpha$, we add the adjustment term given by \eqref{eq:control_adjust} giving us
\vspace{-5pt}
\begin{align*}
    \delta u_{\infty}(x,\bar\alpha,\delta \alpha) = -\sign(g(x,\bar \alpha)) \bigg(\frac{\frac{\partial m (x,\alpha)}{\partial \alpha}}{2\sqrt{m(x,\alpha)}}\bigg)\bigg|_{\alpha = \bar \alpha} \delta\alpha.
\end{align*}

To examine the performance numerically, we consider $g(x,\alpha) = \alpha x^2$ and $m(x,\alpha) = \alpha x^2+ \alpha^2 x^4$, where the nominal value of $\alpha$ is set to $1$. 

The analytical solution is given by
\begin{align*}
    u_{\infty}(x,\alpha) = -\sqrt{\alpha x^2+ \alpha^2 x^4},
\end{align*}
and the variation in the optimal control around the nominal value of $\alpha$ is given by
\begin{align*}
    \delta u_{\infty}(x,\bar\alpha,\delta \alpha) = -\frac{|x|(1+2 \bar \alpha x^2)} {2\sqrt{ \bar\alpha + \bar\alpha^2 x^2}} \delta \alpha.
\end{align*}

For the numerical algorithm using the Galerkin method, the identical set of basis functions used in Example 1 is employed. An initial admissible control input of $u_0 = -x^2$ is chosen. The numerical algorithm is run for 100 iterations, and the optimal controls achieved using analytical and numerical methods are compared shown in Fig. \ref{fig:example2_result} for $\delta \alpha = 0.2$. Similar to the previous example, the NEOC solution matches well with the recalculated solution for both analytical and numerical cases, with the latter
case being visually indistinguishable. The analytical solution reveals that the optimal control law has a non-smooth point at $x=0$. Due to this non-smoothness, when the numerical solution is approximated using polynomial basis functions, significant errors are observed around $x=0$. The difference between the recalculated and NEOC laws is shown in Fig. \ref{fig:example2_error} for different values of parameter perturbations. We observe that even for large perturbations, the NEOC error remains very small.

\begin{figure}[t]
    \centering
    \begin{subfigure}[b]{0.4\textwidth}
        \centering
        \includegraphics[width=\linewidth,trim={0.8cm 0.2cm 0.2cm 1.5cm},clip]{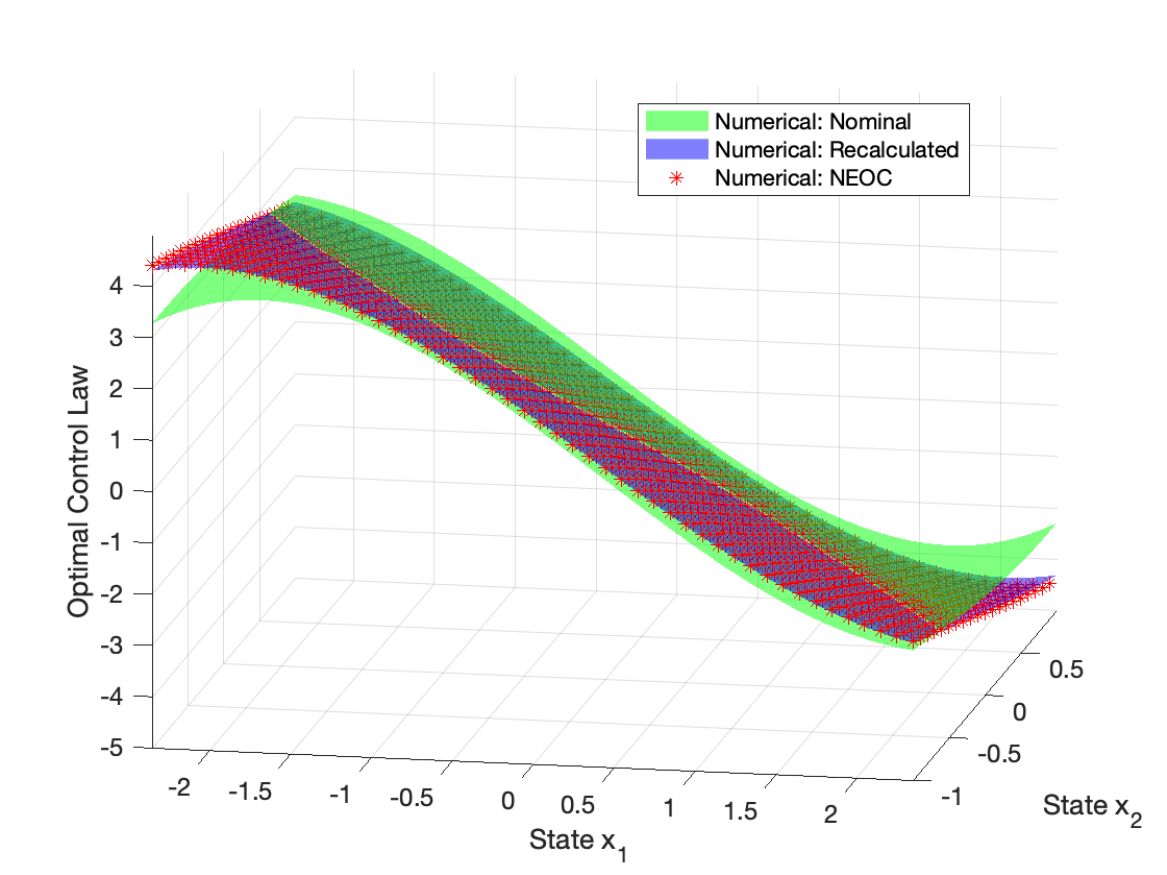}
        \caption{NEOC performance}
        \label{fig:example3_results}
    \end{subfigure}
    \hfill
    \begin{subfigure}[b]{0.4\textwidth}
        \centering
        \includegraphics[width=\textwidth,trim={0.8cm 0.2cm 0.2cm 1.0cm},clip]{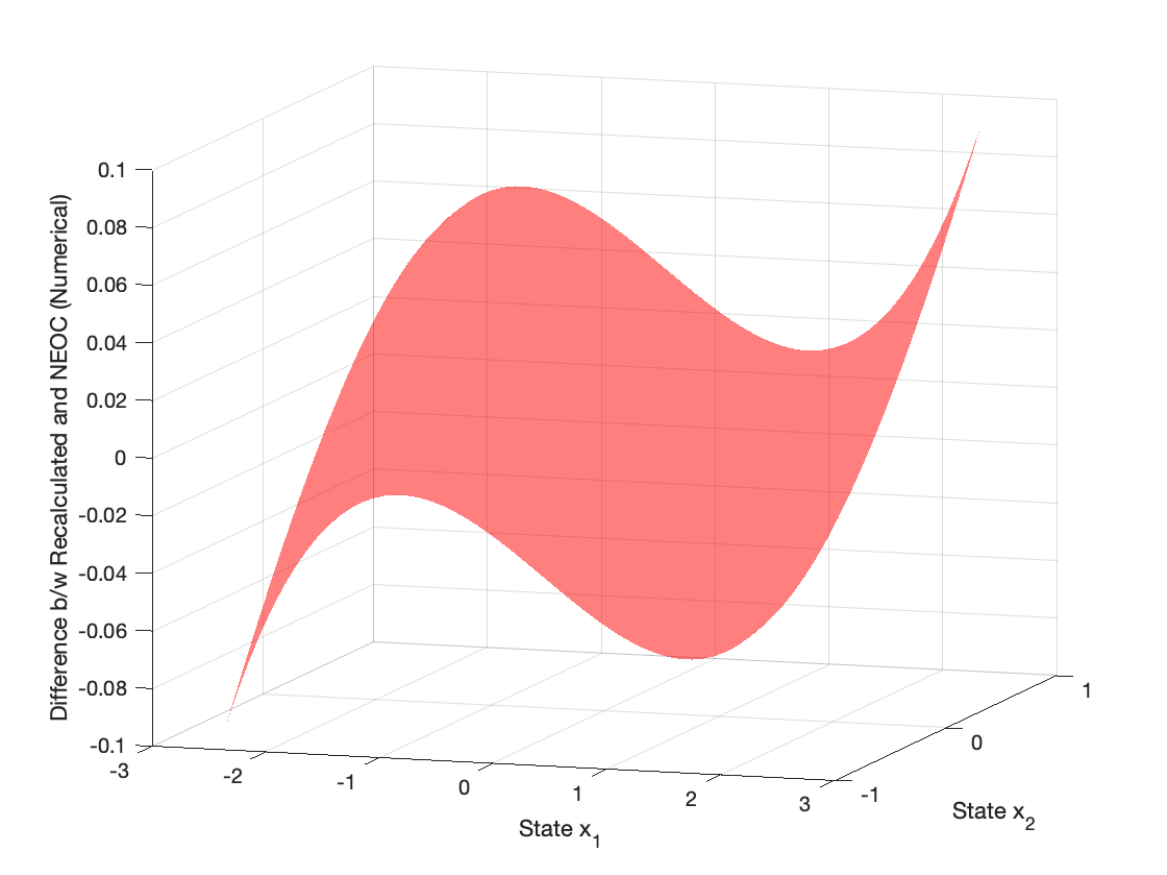}
        \caption{NEOC error}
        \label{fig:example3_error}
    \end{subfigure}
    \caption{Example 3- (a) Comparison of numerical results for nominal, recalculated, and NEOC optimal control law. (b) Difference between the recalculated optimal control law and NEOC for $\delta\alpha = [-0.2 ; 0.2]^{\top}$.}
    \label{fig:example3}
\end{figure}

\begin{example}
Next, we consider a complex dynamics of an inverted pendulum with opposing friction and torque as a general example, which is given by the following equation
\begin{align*}
    m\Ddot{\theta} = mg\sin(\theta)-kl\dot{\theta}+u.
\end{align*}
where $m$ and $l$ represent the mass and length of the rod, $g$ is the gravitational constant, and $k$ is the friction coefficient. By setting $g$ and $k$ to 1 (with no essential loss of generality), the equations of motion can be transformed into a first-order vector system with states $x_1 = \theta$ and $x_2 = \dot{\theta}$, given by
\begin{align*}
    \begin{bmatrix}
    \dot x_1 \\ \dot x_2  
    \end{bmatrix} = \begin{bmatrix}
    x_2 \\ sin(x_1)- \frac{m}{l}x_2  
    \end{bmatrix} + \begin{bmatrix}
    0 \\ 1/m  
    \end{bmatrix}u.  
\end{align*}
We consider two perturbation parameters, denoted by $\alpha = [m ;l]^{\top}$, with their nominal values set as $\bar\alpha = [1 ;1]^{\top}$. The perturbations are defined by $\delta\alpha = [\delta m ; \delta l]^{\top}$, and we analyze the following cost function:
\begin{align*}
V=\lim_{T\to\infty}\int_0^T\big[u^2(x,\alpha)+l^2\|x\|^2]dt.
\end{align*}
\end{example}
The last term denotes the cost associated with the linear position and linear velocity. In this case, obtaining an analytical solution from the equation \eqref{eq:updatephi} is challenging, and hence, we only evaluate the numerical performance by comparing the numerical solution with actual and NEOC values. The basis functions used to approximate the minimum performance index are given by
\vspace{-2pt}
\begin{align*}
    \{\psi_j\} =& \{x_1^2 ,x_1x_2 , x_2^2 , x_1^4,x_1^3x_2, x_1^2 x_2^2, x_1x_2^3, x_2^4 \}.
\end{align*}
The initial admissible control is obtained by canceling the non-linearities in the open-loop system and using an additional stabilizing linear control
\begin{align*}
    u_0(x) = -\sin(x_1) -x_1 - x_2. 
\end{align*}
The surface plots of optimal control laws in Fig. \ref{fig:example3_results} show the performance of the NEOC approach for a perturbation of $\delta \alpha = [-0.2, 0.2]^{\top}$ in the values of the mass and length of the inverted pendulum. The error between the recalculated and NEOC laws is shown in Fig. \ref{fig:example3_error}. It can be observed that the NEOC solution matches very well with the recalculated numerical results, indicating the effectiveness of the NEOC approach in handling perturbations in multiple system parameters.

\section{Homotopic Approach for large perturbation}
\label{sec:homotopic}
As discussed earlier, the error between the NEOC and recalculated optimal control is a function of the magnitude of $\delta\alpha$, the variation in the parameter. Therefore, small perturbations become an essential assumption of our framework due to the first variation condition as large perturbation results in significant errors. Nevertheless, real-world nonlinear problems may necessitate significant perturbations.

In general terms, homotopy is a concept from topology \cite{hatcher2002algebraic} that establishes an equivalence relationship. The fundamental concept involves beginning with a simpler problem (in our case, small perturbations), for which finding a solution is relatively simple. Then, via iterative steps, the solution is transformed into a solution for a more complex problem (involving large perturbations), accomplished by adjusting the homotopy parameter. This idea has been popularly used in the domain of nonlinear differential equations using a homotopy perturbation technique where a homotopy is constructed with an embedding parameter \cite{he1999homotopy,he2003homotopy}. It has also found its application in many optimal control problems \cite{caillau2012differential,bonilla2010convexity,rostalski2011numerical}. In this section, we present a homotopic numerical approach to our NEOC problem that allows us to deal with large perturbations. An approach of this general character for perturbations arising in an open-loop optimal control problem was used in \cite{jiang2015optimal}.
\begin{remark}
    By way of caution, one should note that a change in the optimal control law due to perturbation in an underlying parameter of the optimal control problems \textit{may not always be smooth} (continuous) with respect to the parameter. A simple example to understand this possibility is to consider the case  when a small change in a parameter $\alpha$ makes the system uncontrollable. In scenarios like these, the homotopy method will be ineffective in yielding close-to-optimal solutions.
\end{remark}

Recall that in the numerical approach presented in Section \ref{sec:algorithm}, the NEOC solution is achieved by solving \eqref{eq:u_neoc_algo}, which can correspond to solving NEOC in a single step with the parameter change $\delta \alpha$. If the change is large, then using a  homotopic approach, we decompose this change in parameter into $N$ different steps, each corresponding to a small change. Compared to the single-step case of NEOC which required solving a single linear equation, the homotopic approach requires solving $N$ linear equations. This can be illustrated as  
\begin{align*}
    &\text{Step 1:} \quad u^*(\bar \alpha) \longrightarrow u_N(\bar \alpha + \delta \alpha /N),\\
    &\text{Step 2:} \quad u_N(\bar \alpha+ \delta \alpha /N) \longrightarrow u_N(\bar \alpha + 2\delta \alpha /N),\\
    & \quad\quad\; \vdots \\
    &\text{Step N:} \quad u_N(\bar \alpha+ (N-1)\delta \alpha /N) \longrightarrow u_N(\bar \alpha + \delta \alpha).
\end{align*}

The selection of $N$, representing the number of solution divisions, is a crucial decision. Expanding on this aspect, we present a minimum bound on $N$ in the following corollary derived from Theorem \ref{theorem3}.

\begin{corollary}
    Consider a NEOC problem for which the variation of the optimal control law with $\alpha$ is continuous. Using the same condition as in Theorem \ref{theorem3}, for a given perturbation of $\delta \alpha$, the number of steps required to achieve a desired accuracy of $\epsilon$ is given by
    \begin{align*}
        N \geq \frac{M(\delta\alpha)^2}{2\epsilon},
    \end{align*}
    where $M$ is defined as
\begin{align*}
    M = \max_{x\in \Omega} \max_{\alpha\in [\bar \alpha, \bar \alpha+\delta \alpha]} \left\|\frac{\partial^2 u}{\partial \alpha^2}(x,\alpha)\right\|.
\end{align*}
\end{corollary}
\vspace{5pt}
\noindent{\bf Proof:} With a homotopic approach using $N$ steps, from Theorem 4.3, for any step $i$ we get
\begin{align*}
    \|\delta \alpha_i\| \leq \sqrt{\frac{2\epsilon_i}{M}}.
\end{align*}
By setting the perturbation and accuracy for each step as $\delta \alpha/N$ and $\epsilon/N$ respectively, we get
\begin{align*}
    \frac{\|\delta \alpha\|}{N} \leq \sqrt{\frac{2\epsilon}{NM}},
\end{align*}
resulting in the desired bound. 
$\hfill \blacksquare$

\begin{figure}[htb!]
    \centering

    \includegraphics[width=0.45\textwidth,trim={0.8cm 0.2cm 0.2cm 1.0cm},clip]{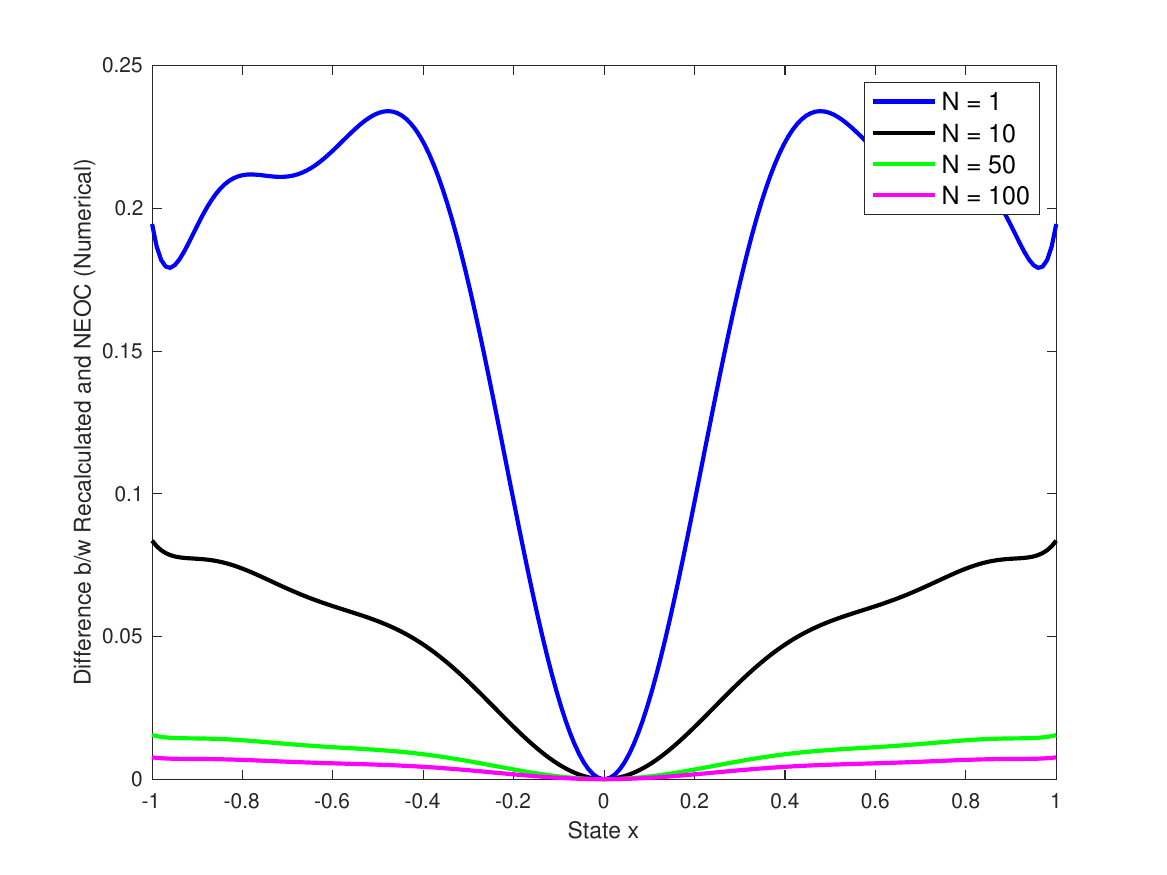}
 
    \caption{Difference between the recalculated optimal control law and NEOC for a different number of steps using the Homotopic approach.}
    \label{fig:homotopy}
\end{figure}
To present the numerical outcomes, we revisit Example 7.3 from Section \ref{sec:examples}. In Fig. \ref{fig:example2_error}, we have depicted the error between the numerical outcomes of the recalculated optimal control law and the NEOC law for different values of perturbation. The highest perturbation considered in the graph is $\delta\alpha = 0.5$, using $\bar \alpha = 1$. Now, we consider a scenario where the desired $\delta\alpha$ amounts to $5$, marking a substantial change of $400\%$. In Fig. \ref{fig:homotopy}, we provide the results obtained through the homotopic approach for varying values of $N$ $(1,10,50,100)$. As anticipated by the theoretical results, we observe a significant reduction in errors as we increase the number of steps, indicating the homotopic approach may be very powerful in practical applications.

\section{Conclusions}
\label{sec:conclusion}
The problem of neighboring extremal optimal control (NEOC) holds substantial relevance across diverse applications. This study aims to address the NEOC challenge in scenarios where the nominal solution involves a closed-loop feedback law, rather than an open-loop control associated with specific initial conditions. Our approach involves a comprehensive analysis of how perturbations in a parameter impact the optimal control law  by studying the variations in the optimal performance index. We derive a closed-loop adjustment term that complements the nominal closed-loop feedback law, facilitating the attainment of neighboring extremal solutions and an overall enhancement in system performance in the presence of small known parameter variations or perturbations. For numerical implementation across general dynamics, we introduce an algorithm that leverages the Galerkin iterative method for solving the Hamilton-Jacobi equations. To illustrate the algorithm's effectiveness and versatility across various scenarios, we provide several examples. Furthermore, we propose a homotopic approach to handle larger perturbations for any given perturbation and desired accuracy. 

Our future efforts will focus on several key areas: 1) the development of efficient online algorithms for solving NEOC in the context of closed-loop control, 2) the creation of a NEOC solution based on closed-loop feedback laws tailored to differential games with performance indices that may not necessarily be quadratic in nature, 3) the use of NEOC to assist in the design of closed-loop nonlinear feedback for linear systems to achieve better-than-linear design tradeoff between rise time and overshoot in a step response, and 4) the implementation of adaptive control strategies, most likely of an iterative identification and feedback redesign variety \cite{albertos2012iterative}, to handle systems in which parameters drift. It will be important both for ensuring that NEOC methods do not give rise to errors that are too large and to avoid the potential for instability in iterative identification and control redesign \cite{anderson2005failures}, to ensure that big parameter changes between successive redesigns do not occur.

%%%%%%%%%%%%%%%%%%%%%%%%%%%%%%%%%%%%%%%%%
\section*{References}
\bibliographystyle{IEEEtran}
\bibliography{sample}
\end{document}